\documentclass[preprint,3p,sort&compress,final,times]{elsarticle}
\usepackage{amsmath}
\usepackage{amssymb}
\usepackage{amsthm}
\usepackage{mathtools}
\usepackage{placeins}
\usepackage{bm}
\usepackage{caption}
\usepackage{subcaption}
\usepackage[usenames,dvipsnames]{color}

\usepackage{cases}

\usepackage{tikz-cd}
\usepackage{soul}

\renewcommand{\d}{\mathrm{d}}

\newcommand \nn {\bm{n}}
\newcommand \uu {\bm{u}}

\newcommand{\tabref}[1]{Table~\ref{#1}} % reference a table
\newcommand{\figref}[1]{Figure~\ref{#1}} % reference a figure
\newcommand{\secref}[1]{Section~\ref{#1}} % reference a section
\newcommand{\exampleref}[1]{Example~\ref{#1}} % reference an example
\newcommand{\transpose}{\intercal}                      % transpose operator
\newcommand{\bb}[1]{\left( #1 \right)}

\newcommand{\bs}[1]{\left[ #1 \right]}
\newcommand{\nl}{\\[1.1ex]}
\newcommand{\conn}{\mathbb{N}_2^\transpose}

\newcommand{\spaces}{representation }
\newcommand{\representations}{representations }

\newcommand{\VJ}[1]{\textcolor{black}{#1}}
\theoremstyle{plain}

\newtheorem{corollary}{Corollary}

\newtheorem{example}{Example}

\newtheorem{lemma}{Lemma}

\journal{}
\begin{document}
\begin{frontmatter}

\title{Use of algebraic dual representations in domain decomposition methods for Darcy flow in 3D domains}
\author[TUD]{V.~Jain\corref{cor}}
\ead{V.Jain@tudelft.nl}
%\author[TUD]{Y.~Zhang}
\author[TUD]{A.~Palha}
\author[TUD]{M.~Gerritsma}

\address[TUD]{Delft University of Technology, Faculty of Aerospace Engineering, P.O. Box 5058, 2600 GB Delft, The Netherlands}
%\address[EUT]{Eindhoven University of Technology, Department of Mechanical Engineering, P.O. Box 513, 5600 MB Eindhoven, The Netherlands}

\cortext[cor]{Corresponding author. Tel. +31 15 2789670.}

\begin{abstract}
In this work we use algebraic dual representations in conjunction with domain decomposition methods for Darcy equations.
We define the broken Sobolev spaces and their finite dimensional counterparts.
In addition, a global trace space is defined that connects the solution between the broken spaces.
Use of dual representations results in a sparse metric free representation of the constraint on divergence of velocity, \VJ{the pressure gradient term} and on the continuity constraint across the sub domains.
%\st{This leads to algebraic systems that have a lower condition number}.
To demonstrate this, we solve two test cases: i) manufactured solution case, and ii) industrial benchmark reservoir modelling problem SPE10.
The results demonstrate that the domain decomposition scheme, although with more unkowns, requires less memory and simulation time as compared to the continuous \VJ{Galerkin} formulation.
\end{abstract}


\end{frontmatter}

\section{Introduction}
Since the early days of computational fluid dynamics there has been an immense increment in computational power.
However, the need for development of fast numerical algorithms has persisted consistently.
To that motive, the class of domain decomposition (DD) methods has played an important role in reducing the run times and memory requirements for numerical simulations.
The general aim for these methods is that the domain of the problem is broken up into a set of many small independent sub domains. 
New degrees of freedom are introduced, the Lagrange multipliers, that impose the appropriate continuity constraints across these domains.
It is then possible to solve for the coupled system of Lagrange multiplier equations only, which is a smaller system compared to the continuous unbroken formulation.
The local solution in the sub domains is obtained independently of each other as a function of the solution of the Lagrange multipliers.
This process reduces the computational burden of the problem in the sense that it is less demanding in terms of required memory and computational time, without compromising on the accuracy of the results.
Over the last decades, many variants of DD methods have been developed, for example, hybrid method~\cite{1977Raviart}, mortar method~\cite{2000Wohlmuth}, FETI method~\cite{1994farhat,2001farhat}, DPG method~\cite{2010Demkowicz,2011Demkowicz,2016Carstensen}, Steklov-Poincar\'e method~\cite{1988Agoshkov,1991Quarteroni}, etc.
%, DtN method  \VJ{add citation},  boundary element method \VJ{add citation}, etc.
For a comprehensive discussion of these methods we refer the reader to the review papers~\cite{2015Cockburn,2016Pietro}.
The two main challenges for DD methods are: i) to find a subset of suitable finite dimensional spaces such that the system matrices are not singular, i.e. they do not produce spurious kinematic modes~\cite{2016Almeida}, especially with respect to the trace spaces of the Lagrange multipliers; and ii) how to efficiently solve the global system of Lagrange multiplier equations which can become large for practical applications and is characterized by high condition numbers.
The primary objective of this paper is to address the first challenge, i.e. to define the framework and extend the use of algebraic dual representations introduced in \cite{2020jain} for DD formulation of Darcy flow.
\VJ{The DD formulation used in this work is based on the hybrid method which is a form of discontinuous Galerking formulation.}
%In particular, we use , that obey the de Rham sequence,  and extend them to the framework of DD methods.
%In \cite{2020jain} the construction of algebraic dual spaces that obey the de Rham sequence is presented.
It is shown that the use of algebraic dual representation results in a sparse metric free representation for the divergence constraint on the velocity, the pressure gradient and the continuity constraint across the sub domain, even for high order spectral element method.
We also see that the continuity constraints are local to the boundary face of adjoining sub domain elements.
This construction is similar to the use of dual spaces in~\cite{2000Wohlmuth}, where they have non-confirming sub domain boundaries and the continuity constraints are also local to the elements on the adjoining sub domain boundaries only.

In this paper we use spectral elements of order $N=1,2,3$ and demonstrate the advantage of using DD formulation with algebraic dual representations using two test cases.
The first test case is a manufactured solution taken from \cite{2012Wheeler} on a randomly deformed smooth domain.
Using this case we first show that using the DD formulation gives the same solution as the continuous formulation and have speed-up in simulation times.
It is observed that on mesh refinement the speed-up in simulation time is larger in the case of high order elements (see \tabref{tab:time_new}).
We also show optimal rates of convergence for varying mesh refinements.
% and elements of order $N = 1,2,3$.
The second test case is a benchmark reservoir modelling test case SPE10 \cite{spe10} that demonstrates the use of this formulation on a practical application.
The purpose of this work is fourfold: i) to demonstrate that the algebraic dual representations can be extended to the framework of DD methods without compromising on the accuracy of results, ii) to show that using DD method we can go to higher levels of mesh refinement with same memory hardware, iii) the use of DD formulation is more efficient in terms of simulation time, and iv) to show the applicability of the method for industry benchmark problems such as the SPE10 case.

This paper is structured as follows:
In \secref{sec:sobolev_spaces} we define the broken Sobolev spaces for DD formulation.
In \secref{sec:finite_spaces} we define the finite dimensional subset of the broken Sobolev spaces.
In \secref{sec:model_problem} we state the Darcy problem and the weak formulation of this problem for DD method.
The algebraic formulation for this problem and the solution steps are also described in this section.
In \secref{sec:test_cases} we present the results for the two test cases, i) the manufactured test case from~\cite{2012Wheeler}, and ii) the benchmark test case SPE10~\cite{spe10}.
We draw conclusions and discuss the scope for future work in \secref{sec:conclusions}.
\section{Broken Sobolev spaces} \label{sec:sobolev_spaces}
Let $\Omega \subset \mathbb{R}^3$ be a bounded domain with Lipschitz boundary $\partial \Omega$.
Let $L^2\bb{\Omega}$ be the space of square integrable functions and $[L^2\bb{\Omega}]^3$ the space of square integrable vector fields in 3D, and $H^1\bb{\Omega}$, $H\bb{\mathrm{div};\Omega}$ be the Hilbert spaces defined as
\begin{equation*}
\begin{array}{ll}
H^1\bb{\Omega} & := \left\lbrace p \in L^2\bb{\Omega} : \mathrm{grad}\ p \in [L^2\bb{\Omega}]^3 \right\rbrace \nl
H\bb{\mathrm{div}; \Omega} & :=\left\lbrace \bm{v} \in [L^2\bb{\Omega}]^3 : \mathrm{div}\ \bm{v} \in L^2\bb{\Omega} \right\rbrace
\end{array} \;.
\end{equation*}
The trace spaces of $H^1\bb{\Omega}$ is defined as
\begin{equation*}
\begin{array}{ll}
H^{1/2}\bb{\partial \Omega} & := \left\lbrace \lambda \in L^2\bb{\partial \Omega} : \exists \ p \in H^1\bb{\Omega}\quad \mathrm{s.t.}\quad \lambda = p|_{\partial \Omega}  \right\rbrace \nl
%H^{-1/2}\bb{\partial \Omega} & := \left\lbrace \mu \in H^{-1}\bb{\partial \Omega} : \exists \ \bm{v} \in H\bb{\mathrm{div};\Omega} \quad \mathrm{s.t.} \quad \mu = \bm{v}|_{\partial \Omega} \cdot \bm{n} \right\rbrace
\end{array} \;,
\end{equation*}
and we denote by $H^{-1/2}\bb{\partial \Omega}$ its dual space.

Let $\Omega$ be broken into $T$ non-overlapping open sub domains $   \mathcal{M} _i $ with Lipschitz boundary $\partial \mathcal{M}_i$, $i = 1 \;, \hdots \;, T$, such that
\begin{equation} \label{eq:break_domain}
\Omega = \bigcup _{i=1, \hdots \;, T} \mathcal{M}_i \qquad \mbox{with,} \qquad \mathcal{M}_i \subset \Omega \qquad \mbox{and} \qquad \overset{\circ}{\mathcal{M}}\ _i \bigcap \overset{\circ}{\mathcal{M}}\ _j = \emptyset \quad \mbox{for} \quad i \neq j \;.
\end{equation}
Let $\Omega _T$ be the set of sub domains, and $\partial \Omega _T$ be the set of boundaries of these sub domains defined as
\[ \Omega _T = \left\lbrace \mathcal{M}_i \right\rbrace _{i=1, \hdots \;, T} \qquad \qquad \partial \Omega _T =  \left\lbrace \partial \mathcal{M}_i \right\rbrace _{i=1, \hdots \;, T} \;. \] 
%We will denote the interface between two domains $\mathcal{M}_i$ and $\mathcal{M}_j$ as $\Gamma _{ij}$, such that 
%\[ \Gamma _{ij} = \partial \mathcal{M} _i \bigcap \partial \mathcal{M}_j \qquad \mbox{for} \qquad 1 \leq i < j \leq T \;, \]
%and the collection of interfaces as
%\[\VJ{ \Gamma _T = {\left\lbrace \Gamma _{ij} \right\rbrace}_{1 \leq i < j \leq T} } \;. \]
We define the broken Sobolev spaces for the set of sub domains $\Omega _T$ as
\begin{eqnarray}
L^2\bb{\Omega _T} := &   & \prod _{\mathcal{M} \in \Omega _T} L^2\bb{\mathcal{M}} \nonumber \nl
%H^1\bb{\Omega _T} := & \left\lbrace p|_\mathcal{M} \in H^1\bb{\mathcal{M}}, \mathcal{M} \in \Omega _T \right\rbrace & = \prod _ {\mathcal{M} \in \Omega _T} H^1\bb{\mathcal{M}} \nonumber \nl
H\bb{\mathrm{div};\Omega _T} := &  & \prod _{\mathcal{M} \in \Omega _T} H\bb{\mathrm{div \;; \mathcal{M}}} \;. \nonumber \nl
H^{-1/2} \bb{\partial \Omega _T} := & & \prod _{\Gamma \in \partial \Omega _T} H^{-1/2} \bb{\Gamma}  \nonumber
\end{eqnarray}
Let the set of boundary faces, $\partial \Omega _e$ be defined as
\[ \partial \Omega _e = \left\lbrace \gamma _{ij} = \partial \overline{ \mathcal{M}}_i \bigcap \partial \overline{\mathcal{M}}_j  \qquad \mbox{for} \qquad i<j \right\rbrace \;. \]
\[ \partial \Omega _\Gamma = \left\lbrace \gamma_{i} =\partial \mathcal{M}_i \bigcap \partial \Omega \qquad \mbox{for} \qquad i = 1, \hdots , T  \right\rbrace \;. \]
%Let $e$ be the total number of unique boundary faces in the set $\partial \Omega _T$, and $\partial \Omega _e$ be the set of all such unique boundary faces.
Then we define the trace space
\[ H^{1/2} \bb{\partial \Omega _e} = \left\lbrace \lambda |_{\gamma} \in L^2\bb{\gamma}, \gamma \in \partial \Omega _e \right\rbrace = \prod _{\Gamma \in \partial \Omega _e} H^{1/2}\bb{\Gamma} \;.  \]
\subsection*{Notation}
We will denote the $L^2$-inner product by $\bb{\cdot \;, \cdot}$.
A set of coordinates is denoted by $\bm{x}$.
We denote by $N$ the highest polynomial degree of basis functions used in an element is $N-1$.

We use two set of finite element spaces, the primal representation and the dual representation from \cite{2020jain}.
The primal representation of finite element space, basis functions, and degrees of freedom are dentoed by $X\bb{\cdot}$, $\Psi\bb{\cdot}$, $\mathcal{N}$ respectively.
The dual representation of space $X\bb{\cdot}$, its associated basis functions and expansion coefficients are denoted with a tilde as $\widetilde{X}\bb{\cdot}$, $\widetilde{\Psi}\bb{\cdot}$ and $\widetilde{\mathcal{N}}\bb{\cdot}$, respectively.
The basis functions are always represented as row vectors 
\[ \Psi ^k \bb{\bm{x}} = \bb{\bm{\epsilon} ^{(k)}_1\bb{\bm{x}} \quad \bm{\epsilon} ^{(k)}_2\bb{\bm{x}} \quad \hdots \quad \bm{\epsilon} ^{(k)}_{d_k}\bb{\bm{x}}} \;, \]
where $d_k$ is the dimension of polynomial vector space of $X\bb{\cdot}$.
The expansion coefficients are always represented as column vectors
\[ \mathcal{N}^k\bb{\bm{u}} = \bb{\begin{array}{c}
\mathcal{N}_1^k\bb{\bm{u}} \nl
\mathcal{N}_2^k\bb{\bm{u}} \nl
\vdots \nl
\mathcal{N}_{d_k}^k\bb{\bm{u}}
\end{array}} \]
We will Lemma 2 from \cite{2020jain}, which states
\begin{lemma}
If $\Psi^k\bb{\bm{x}}$ and $\widetilde{\Psi}^k\bb{\bm{x}}$ are basis functions from primal and dual representations respectively, then these bases are bi-orthogonal with respect to each other
\begin{equation} \label{eq:primal_dual_int}
\int _\Omega \widetilde{\Psi}^{3-k} \bb{\bm{x}} \Psi ^k \bb{\bm{x}} \mathrm{d}\Omega = \mathbb{I} \;,
\end{equation}
where $\mathbb{I}$ is the identity matrix of dimension $d_k$.
\end{lemma}
\begin{corollary} \label{cor:prop1}
The inner product between variables of primal and dual representation is vector dot product of expansion coefficients.
\begin{proof}
Let $p,q \in X\bb{\Omega} \times \widetilde{X}\bb{\Omega}$, the the inner product is given by
\[ \int_\Omega p\ q\ \mathrm{d}\Omega =  \widetilde{\mathcal{N}}^{3-k}\bb{q}^\transpose \bb{\int_\Omega \widetilde{\Psi} ^{3-k} \bb{\bm{x}}^\transpose \Psi ^k \bb{\bm{x}} \mathrm{d}\Omega} \mathcal{N}^k\bb{p} \stackrel{\eqref{eq:primal_dual_int}}{=} \widetilde{\mathcal{N}}^{3-k}\bb{q}^\transpose \mathcal{N}^k\bb{p} \;. \]
\end{proof}
\end{corollary}
In general, if not explicitly mentioned otherwise, we use Gauss-Lobatto-Legendre points for numerical integration.
%\VJ{
%The expansion coefficients for the finite dimensional trace spaces are denoted by $\mathcal{B}\bb{\cdot}$ instead of $\mathcal{N}\bb{\cdot}$, and the corresponding dual expansion coefficients are denoted by $\dual{\mathcal{B}}\bb{\cdot}$.
%}
\section{Finite dimensional spaces} \label{sec:finite_spaces}
In this section we will define the finite dimensional spaces for hexahedral elements in 3D domains.
The domain $\Omega$ is discretized into a mesh that consists of points, edges, surfaces and volumes.
For $N>1$ each element also consists of a GLL mesh.
We will first introduce the finite dimensional spaces for the domain $\Omega$ as defined in \cite{2020jain}.
The basis functions for domain $\Omega$ are obtained from mapping of basis functions on a reference domain $\widehat{\Omega} = [-1,1]^3$.
For further details on construction we refer the reader to \cite[\S 4.5]{2020jain}.
Here we directly use the mapped basis functions and the relevant properties of these spaces.
We will use the primal representations to define finite element spaces: i) $D\bb{\Omega} \subset H\bb{\mathrm{div};\Omega}$, ii) $S\bb{\Omega} \subset L^2\bb{\Omega}$, iii) $D_b\bb{\partial \Omega} \subset H^{-1/2}\bb{\partial \Omega}$, and  dual representations to define the finite element spaces: i) $\widetilde{S}\bb{\Omega} \subset L^2\bb{\Omega}$, ii) $\widetilde{D}_b\bb{\partial \Omega} \subset H^{1/2}\bb{\partial \Omega}$.

These spaces will then be used to define the broken finite dimensional spaces that will be used for DD formulation.
\subsection{Primal representations}
%\VJ{Let $N_f$ be the total number of boundary faces and $N_v$ be the total number of volumes in the discretized domain $\Omega$.}
\subsubsection{Finite element space $D\bb{\Omega} \subset H\bb{\mathrm{div};\Omega}$}
Let $\bm{u}$ be the flux component of a vector field. 
For an element $\bm{u} \in D\bb{\Omega}$, let $\Psi^2\bb{\bm{x}}$, $\mathcal{N}^2\bb{\bm{u}}$ form the row vector of basis functions and the column vector of expansion coefficients, given as
\begin{equation} \label{eq:basis_surface}
\Psi ^2 \bb{\bm{x}} = \bb{\bm{\epsilon} ^{(2)}_1\bb{\bm{x}} \quad \bm{\epsilon} ^{(2)}_2\bb{\bm{x}} \quad \hdots \quad \bm{\epsilon} ^{(2)}_{d_2}\bb{\bm{x}}} \;, \qquad \qquad \mathcal{N}^2\bb{\bm{u}} = \bb{\begin{array}{c}
\mathcal{N}_1^2\bb{\bm{u}} \nl
\mathcal{N}_2^2\bb{\bm{u}} \nl
\vdots \nl
\mathcal{N}_{d_2}^2\bb{\bm{u}}
\end{array}} \;, \qquad \qquad \mbox{respectively.}
\end{equation}
Here, $\mathcal{N}^2_i \bb{\bm{u}} = \int _{f_i} \bm{u}\cdot \bm{n}\ \mathrm{dA}$ denotes the net flux $\bm{u}$ through the face $f_i$.

Then we can represent $\bm{u}$ as
\begin{equation*}
\bm{u}\bb{\bm{x}} = \Psi	 ^2 \bb{\bm{x}} \mathcal{N}^2\bb{\bm{u}} \;.
\end{equation*}
Using this notation, for two elements $\bm{u},\bm{v} \in D\bb{\Omega}$ the $L^2$-inner product is given as
\begin{equation*} \label{eq:inner_Dc}
\bb{\bm{u} \;, \bm{v}} = \int _\Omega \bm{u}^\transpose \bm{v} \mathrm{d}\Omega = {\mathcal{N}^2\bb{\bm{u}}}^\transpose \bb{\int _\Omega {\Psi ^2 \bb{\bm{x}}}^\transpose \Psi ^2 \bb{\bm{x}} \mathrm{d} \Omega} \ \mathcal{N}^2\bb{\bm{v}} = {\mathcal{N}^2\bb{\bm{u}}}^\transpose \mathbb{M}^{(2)} \mathcal{N}^2\bb{\bm{v}} \;,
\end{equation*}
where, $\mathbb{M}^{(2)} := \int _\Omega {\Psi ^2 \bb{\bm{x}}}^\transpose \Psi ^2 \bb{\bm{x}} \mathrm{d} \Omega$ is the mass matrix associated with the basis functions $\Psi ^2 \bb{\bm{x}}$.

If there is a symmetric positive definite permeability tensor $\mathbb{K}$, then the weighted inner product is given as
\begin{equation*}
\bb{\bm{u} \;, \bm{v}}_{\mathbb{K}^{-1}} = \int _\Omega \bm{u}^\transpose \mathbb{K}^{-1} \bm{v} \mathrm{d}\Omega = {\mathcal{N}^2\bb{\bm{u}}}^\transpose \bb{\int _\Omega {\Psi ^2 \bb{\bm{x}}}^\transpose \mathbb{K}^{-1}\bb{\bm{x}} \Psi ^2 \bb{\bm{x}} \mathrm{d} \Omega} \ \mathcal{N}^2\bb{\bm{v}} = {\mathcal{N}^2\bb{\bm{u}}}^\transpose \mathbb{M}^{(2)}_{\mathbb{K}^{-1}} \mathcal{N}^2\bb{\bm{v}} \;,
\end{equation*}
where $\mathbb{M}^{(2)}_{\mathbb{K}^{-1}} := \int _\Omega {\Psi ^2 \bb{\bm{x}}}^\transpose \mathbb{K}^{-1}\bb{\bm{x}} \Psi ^2 \bb{\bm{x}} \mathrm{d} \Omega$ is the mass matrix of the weighted inner product.
\subsubsection{Finite dimensional space $S\bb{\Omega} \subset L^2\bb{\Omega}$}
Let $p$ be any element of the space $S\bb{\Omega}$ and $N_v$ be the total number of volumes in the discretized domain $\Omega$.
If $\Psi	^3 \bb{\bm{x}}$ and $\mathcal{N}^3\bb{p}$ form the row vector of basis functions and the column vector of expansion coefficients, given as
\begin{equation*}
\Psi ^3 \bb{\bm{x}} = \bb{\epsilon ^{(3)}_1\bb{\bm{x}} \quad \epsilon ^{(3)}_2\bb{\bm{x}} \quad \hdots \quad \epsilon ^{(3)}_{d_3}\bb{\bm{x}}} \;, \qquad \qquad \mathcal{N}^3\bb{p} = \bb{\begin{array}{c}
\mathcal{N}_1^3\bb{p} \nl
\mathcal{N}_2^3\bb{p} \nl
\vdots \nl
\mathcal{N}_{d_3}^3\bb{p}
\end{array}} \;, \qquad \qquad \mbox{respectively,}
\end{equation*}
%where $\epsilon^{(3)}_1 \;, \epsilon^{(3)}_2 \;, \hdots $ are the basis functions and $\mathcal{N}^3_1 \bb{p} \;, \mathcal{N}^3_2 \bb{p} \;, \hdots$ are the expansion coefficients, 
then we can represent $p$ as
\begin{equation*}
p \bb{\bm{x}} = \Psi^3 \bb{\bm{x}} \mathcal{N}^3 \bb{p} \;,
\end{equation*}
where $\mathcal{N}^3_i\bb{p} = \int_{V_i} p\ \mathrm{d}V$ denotes the integral of $p$ over mesh volume $V_i$. 

Using this notation, the $L^2$-inner product for two elements $p,q \in S\bb{\Omega}$ is given as
\begin{equation} \label{eq:inner_S}
\bb{p \;, q} = \int _\Omega p^\transpose q\ \mathrm{d}\Omega = \mathcal{N}^3 \bb{p}^\transpose \bb{\int _\Omega {\Psi ^3 \bb{\bm{x}}}^\transpose \Psi ^3\bb{\bm{x}} \mathrm{d}\Omega} \ \mathcal{N}^3 \bb{q} = \mathcal{N}^3 \bb{p}^\transpose \mathbb{M}^{(3)} \ \mathcal{N}^3 \bb{q} \;,
\end{equation}
where, $\mathbb{M}^{(3)} := \int _\Omega {\Psi ^3 \bb{\bm{x}}}^\transpose \Psi ^3 \bb{\bm{x}} \mathrm{d}\Omega$ is the mass matrix associated to the basis functions $\Psi ^3 \bb{\bm{x}}$.
\subsubsection{$D_b\bb{\partial \Omega} \subset H^{-1/2}\bb{\partial \Omega}$}
The trace space $D_b\bb{\partial \Omega}$ is defined as the restriction of vector fields in $D\bb{\Omega}$ to the domain boundary $\partial \Omega$.
If $\mathbb{N}_2$ is the discrete representation of the inclusion map, that maps degrees of freedom defined on the boundary to the global degrees of freedom of the element, and $\mathbb{N}_2^\transpose$ is the discrete representation of the trace operator thaat restricts global degrees of freedom to the boundary, then the basis functions on the boundary and degrees of freedom on the boundary are given by
\begin{equation} \label{eq:bdof}
\Psi ^2_b \bb{\bm{x}} := \Psi ^2\bb{\bm{x}} \mathbb{N}_2 \qquad\mbox{and} \qquad \mathcal{B}^2\bb{\bm{u} \cdot \bm{n}} := \mathbb{N}_2^\transpose \mathcal{N}^2\bb{\bm{u}} \;.
\end{equation}
%In \eqref{eq:basis_surface} if we restrict $\bm{x} \in \partial \Omega$ then most of the basis functions are identically zero and we can remove them from the list.
%Likewise we can remove the corresponding degrees of freedom.
%We will label the remaining basis functions by the subscript $b$ and degrees of freedom by $\mathcal{B}$.
%The construction of trace space $D_b\bb{\partial \Omega}$ is the same as that of 2D finite dimensional space $S\bb{\Omega}$ described in \cite[\S 3.2]{2020jain}.

If $\mathrm{tr}\bm{u} \in D_b\bb{\partial \Omega}$ represents the restriction of flux component of vector field to the domain boundary, then the solution on the boundary can be represented by
\[ \mathrm{tr}\ \bm{u}\bb{\bm{x}} = \Psi ^2_b\bb{\bm{x}} \mathcal{B}^2 \bb{\bm{u}} \;. \]
%For unit outward normal vector fields $\bm{u}\cdot \bm{n} \in D_b\bb{\partial \Omega}$, let $\Psi ^2_b \bb{\bm{x}}$ be the row vector of basis functions and $\mathcal{B}^2\bb{\bm{u} \cdot \bm{n}}$ the column vector of expansion coefficients, 
%where we use $\mathcal{B}$ instead of $\mathcal{N}$ to denote that the expansion coefficients are on the boundary, 
where tr is the trace operator.

%In \eqref{eq:bdof} $\mathbb{N}_2$ is the inclusion matrix that maps the expansion coefficients from the trace space $D_b\bb{\partial \Omega}$ to the expansion coefficients of the space $D\bb{\Omega}$ and \VJ{$\mathbb{N}_2^\transpose$ is the trace matrix that restricts expansion coefficients of $D\bb{\Omega}$ to the space $D_b\bb{\partial \Omega}$}.
For construction of inclusion matrix, see \cite[Ex. 3]{2020jain}.
The inclusion matrix is a sparse metric-free matrix, that consists of $+1,-1$ and $0$ entries only, and is independent of shape and size of elements as long as the topology, numbering of the degrees of freedom and the orientation of the mesh remains the same.
\subsubsection{The divergence operator $\mathbb{E}^{3,2}$}
For any element $\bm{u} \in D\bb{\Omega}$, the divergence operation on $\bm{u}$ is defined as (see \cite[\S 4.4]{2020jain}): $\mathrm{div} : D\bb{\Omega} \longrightarrow S\bb{\Omega}$, such that
\begin{equation} \label{eq:div1}
\mathrm{div} \ \bm{u} = \Psi	^3\bb{\bm{x}} \mathbb{E}^{3,2} \mathcal{N}^2\bb{\bm{u}} \;,
\end{equation}
where $\mathbb{E}^{3,2}$ is the discrete representation of the divergence operator that acts on the expansion coefficients of $\bm{u}$.
The divergence operation changes the degrees of freedom and the basis functions to those of space $S(\Omega)$.
\begin{example} \label{ex:topology}
Divergence operator for 2D mesh
\begin{figure}
\centering
\includegraphics[scale=0.4]{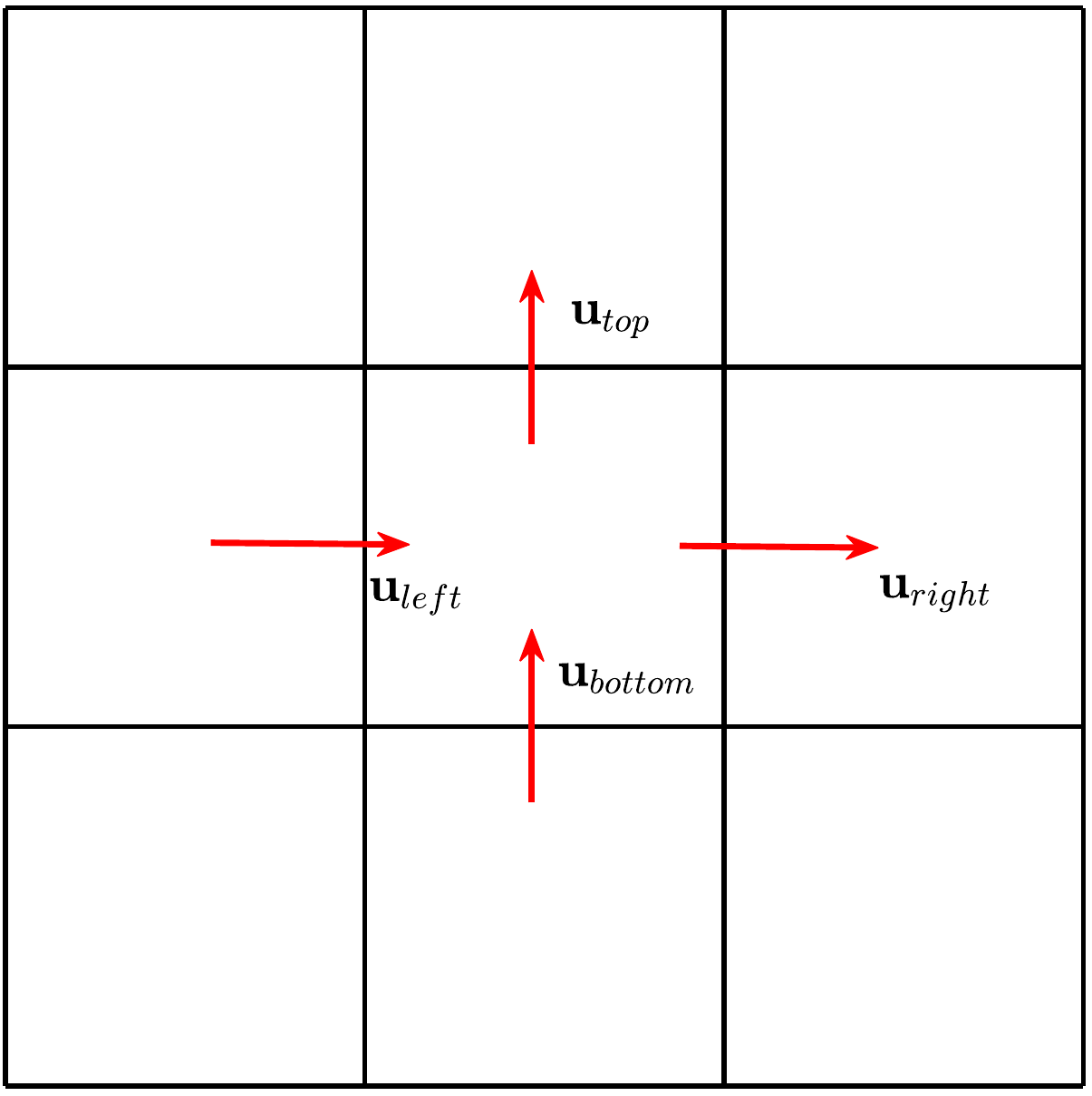} \qquad
\includegraphics[scale=0.5]{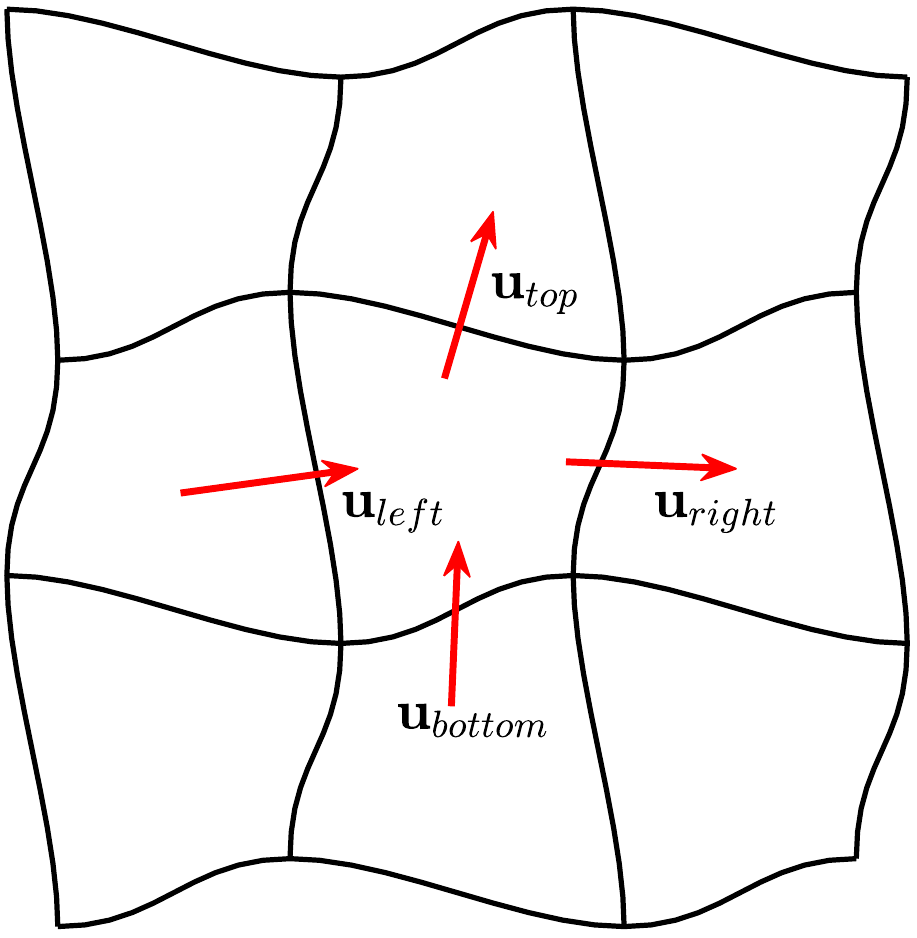}
\caption{Left: A square domain divided into $3 \times 3$ elements. Right: A deformed square domain divided into $3 \times 3$ elements.}
\label{fig:divergence}
\end{figure}

In \figref{fig:divergence} on the left plot we show a square domain divided into $3 \times 3$ elements.
The expansion coefficients of $\bm{u} \in D\bb{\Omega}$, are defined \VJ{across} the edges in the mesh.
The divergence operation on any element $K$, is then defined as
\begin{equation} \label{eq:stokes}
\int _K \mathrm{div} \ \bm{u} \ \mathrm{d}K = \int _{\partial K} \bm{u} \cdot \bm{n} \ \mathrm{d}\gamma = -\mathcal{N}^1\bb{\bm{u}} _{left} + \mathcal{N}^1\bb{\bm{u}} _{right} -\mathcal{N}^1\bb{\bm{u}} _{bottom} + \mathcal{N}^1\bb{\bm{u}} _{top} \;.
\end{equation}
If we assemble \eqref{eq:stokes} for all the nine elements, with appropriate numbering, we get the discrete divergence operator as
\begin{equation*}
%\label{eq:divergence_operator}
\mathbb{E}^{2,1} = 
\resizebox{0.93\hsize}{!}{$
\left[ \begin{array}{cccccccccccccccccccccccc}
-1 & 0 & 0 & 1 & 0 & 0 & 0 & 0 & 0 & 0 & 0 & 0 & -1 & 0 & 0 & 1 & 0 & 0 & 0 & 0 & 0 & 0 & 0 & 0 \\
0  & 0 & 0 & -1 & 0 & 0 & 1 & 0 & 0 & 0 & 0 & 0 & 0 & -1 & 0 & 0 & 1 & 0 & 0 & 0 & 0 & 0 & 0 & 0 \\
0  & 0 & 0 & 0  & 0 & 0 & -1 & 0 & 0 & 1 & 0 & 0 & 0 & 0 & -1 & 0 & 0 & 1 & 0 & 0 & 0 & 0 & 0 & 0 \\
0 & -1 & 0 & 0 & 1 & 0 & 0 & 0 & 0 & 0 & 0 & 0 & 0 & 0 & 0 & -1 & 0 & 0 & 1 & 0 & 0 & 0 & 0 & 0 \\
0 & 0  & 0 & 0 & -1 & 0 & 0 & 1 & 0 & 0 & 0 & 0 & 0 & 0 & 0 & 0 & -1 & 0 & 0 & 1 & 0 & 0 & 0 & 0 \\
0 & 0 & 0 & 0 & 0 & 0 & 0 & -1 & 0 & 0 & 1 & 0 & 0 & 0 & 0 & 0 & 0 & -1 & 0 & 0 & 1 & 0 & 0 & 0 \\
0 & 0 & -1 & 0 & 0 & 1 & 0 & 0 & 0 & 0 & 0 & 0 & 0 & 0 & 0 & 0 & 0 & 0 & -1 & 0 & 0 & 1 & 0 & 0 \\
0 & 0 & 0 & 0 & 0 & -1 & 0 & 0 & 1 & 0 & 0 & 0 & 0 & 0 & 0 & 0 & 0 & 0 & 0 & -1 & 0 & 0 & 1 & 0 \\
0 & 0 & 0 & 0 & 0 & 0 & 0 & 0 & -1 & 0 & 0 & 1 & 0 & 0 & 0 & 0 & 0 & 0 & 0 & 0 & -1 & 0 & 0 & 1
\end{array} \right] \;.
%\bs{\begin{array}{c}
%\overline{\bm{u}}_1 \nl
%\overline{\bm{u}}_2 \nl
%\vdots \nl
%\overline{\bm{u}}_{23} \nl
%\overline{\bm{u}}_{24}
%\end{array}
%} \;.
$}
\end{equation*}
It is a sparse metric free matrix that consists of $+1,-1,0$ entries only.
If the domain is deformed, for example see the right plot of \figref{fig:divergence}, but the connection between the nodes, edges, surfaces, and volumes, remains the same, \VJ{relation} \eqref{eq:stokes} remains the same and consequently the matrix $\mathbb{E}^{2,1}$ remains unchanged.
\end{example}
\VJ{\subsection{Algebraic dual \representations}
In this section we will introduce the finite dimensional sub spaces for $L^2\bb{\Omega}$ and $H^{1/2}\bb{\partial \Omega}$ using the algebraic dual \representations as defined in \cite{2020jain}.}
\subsubsection{Finite dimensional space $\widetilde{S}\bb{\Omega} \subset L^2\bb{\Omega}$}
For the finite dimensional space $S\bb{\Omega}$ let $\widetilde{S}\bb{\Omega}$ be the corresponding algebraic dual \spaces.
For any element $\VJ{q} \in \widetilde{S}\bb{\Omega}$, if $\widetilde{\Psi}^0 \bb{\bm{x}}$ and $\widetilde{\mathcal{N}}^0 \bb{\widetilde{q}}$ are the associated basis functions and the expansion coefficients, then we can represent $q$ as
\[\VJ{q}\bb{\bm{x}} = \widetilde{\Psi}^0\bb{\bm{x}} \widetilde{\mathcal{N}}^0 \bb{\VJ{q}} \;, \]
where
\begin{equation*}
\widetilde{\Psi}^0 \bb{\bm{x}} = \Psi^3 \bb{\bm{x}} \bb{\mathbb{M}^{(3)}}^{-1} \;, \qquad \text{and} \;, \qquad \widetilde{\mathcal{N}}^0 \bb{\VJ{q}} = \mathbb{M}^{(3)} \mathcal{N}^3\bb{q} \;.
\end{equation*}
Using Corollary \ref{cor:prop1}, the $L^2$-inner product between the elements $p,q \in S\bb{\Omega} \times \widetilde{S}\bb{\Omega}$ is then given by
\begin{equation} \label{eq:duality_S}
\bb{q \;, p} = \int _\Omega q^\transpose p \ \mathrm{d}\Omega = {\widetilde{\mathcal{N}}^0\bb{\VJ{q}} }^\transpose \mathcal{N}^3\bb{p} \;.
\end{equation}
We see that the $L^2$-inner product in \eqref{eq:inner_S} requires the evaluation of mass matrix, whereas the \VJ{inner product} in \eqref{eq:duality_S} requires only the vector product between the expansion coefficients which makes the discrete system more sparse and easier to set up.
\subsubsection{Finite dimensional trace space $\widetilde{D}_b\bb{\partial \Omega} \subset H^{1/2}\bb{\partial \Omega}$}
It is known that the Sobolev spaces $H^{1/2}\bb{\partial\Omega}$ and $H^{-1/2}\bb{\partial\Omega}$ are dual to each other.
To replicate this duality also in the discrete setting, we choose $\widetilde{D}_b\bb{\partial \Omega}$ as finite dimensional sub space of $H^{1/2}\bb{\partial\Omega}$ which is the algebraic dual \spaces of $D_b\bb{\partial \Omega} \subset H^{-1/2}\bb{\partial \Omega}$.

For an element $\lambda \in \widetilde{D}_b\bb{\partial \Omega}$ if $\widetilde{\Psi}_b^0 \bb{\bm{x}}$ and $\widetilde{\mathcal{B}}^0\bb{\lambda}$ are the basis functions and the expansion coefficients then, we can represent $\lambda$ as
\begin{equation}
\lambda \bb{\bm{x}} = \widetilde{\Psi}^0_b \bb{\bm{x}} \widetilde{\mathcal{B}}^0\bb{\lambda} \;.
%\qquad \mbox{where,} \qquad \widetilde{\mathcal{B}}^0\bb{\lambda} \VJ{:=} \int _{\partial \Omega} \Psi^2_b \bb{\bm{x}} \lambda \bb{\bm{x}} \mathrm{d}\Gamma \;,
\end{equation}
where
\begin{equation}
\widetilde{\Psi}^0_b \bb{\bm{x}} = \Psi ^2_b \bb{\bm{x}}  \bb{\mathbb{M}_b^{(2)}}^{-1} \qquad \widetilde{\mathcal{B}}^0\bb{\lambda} = \mathbb{M}^{(2)}_b \mathcal{B}^2 \bb{\lambda} \;, \qquad \mbox{with,} \qquad \mathbb{M}_b^{(2)} = \int \Psi ^2_b \bb{\bm{x}}^\transpose \Psi ^2_b \bb{\bm{x}} \mathrm{d}\Gamma \;.
\end{equation}
The inner product between the elements $\lambda \;, \mathrm{tr}\ \bm{u} \in \widetilde{D}_b\bb{\partial \Omega} \times D_b\bb{\partial \Omega}$, is given by
\begin{equation} \label{eq:trace_duality}
\bb{\lambda \;, \mathrm{tr}\ \bm{u}} = \widetilde{\mathcal{B}}^0\bb{\lambda}^\transpose \bb{\int _\Omega \widetilde{\Psi}^0_b \bb{\bm{x}}^\transpose \Psi ^2_b \bb{\bm{x}} \mathrm{d}\Omega } \mathcal{B}^2\bb{\mathrm{tr}\ \bm{u}} = \widetilde{\mathcal{B}}^0\bb{\lambda}^\transpose \mathcal{B}^2\bb{\mathrm{tr}\ \bm{u}} \;.
\end{equation}
\subsection{The gradient of dual \representations}
For a scalar field $p \in \widetilde{S}\bb{\Omega},$ and $\hat{p} \in \widetilde{D}_b\bb{\partial \Omega}$ its values on the domain boundary, the gradient operation for dual \representations is defined as \cite[Def. 19]{2020jain}, $\widetilde{\mathrm{grad}} : \widetilde{S}\bb{\Omega} \times \widetilde{D}_b\bb{\partial \Omega} \rightarrow \widetilde{D}\bb{\Omega}$, such that 
%for $p \;, \lambda \in \widetilde{S}\bb{\Omega} \times \widetilde{D}_b\bb{\partial \Omega}$
\begin{equation}
\int _\Omega \widetilde{\mathrm{grad}} \bb{p \;, \hat{p}} \bm{u} \ \mathrm{d}\Omega = - \int _\Omega p \bb{\mathrm{div}\ \bm{u}} \ \mathrm{d} \Omega + \int _{\partial \Omega} \hat{p} \bb{\bm{u} \cdot \bm{n}} \ \mathrm{d} \Gamma \qquad \qquad \forall \bm{u} \in D\bb{\Omega} \;.
\end{equation}
The expansion coefficients of $\widetilde{\mathrm{grad}}\ \bb{p, \hat{p}}$ are given as
\begin{equation}
\widetilde{\mathcal{N}}^1 \bb{\widetilde{\mathrm{grad}} \bb{{p}\;,\hat{p} }} = - {\mathbb{E}^{3,2}}^\transpose \widetilde{\mathcal{N}}^0\bb{p} + \mathbb{N}_2 \widetilde{\mathcal{B}}^0 \bb{\hat{p}} \;,
\end{equation}
expanded in the basis $\widetilde{\Psi}^1\bb{\bm{x}} = \mathbb{M}^{(2)}\Psi ^2\bb{\bm{x}}$.
\subsection{Broken finite dimensional spaces}
In this section we define the finite dimensional subset of Broken Sobolev spaces.
\subsubsection{Finite dimensional space $D\bb{\Omega _T} \subset H\bb{\mathrm{div};\Omega _T}$}
For the set of domains $\Omega _T$ we define the finite element space $D\bb{\Omega _T} \subset H\bb{\mathrm{div};\Omega _T}$ as
%\begin{equation}
%D\bb{\Omega _T} := \left\lbrace \bm{u} |_{\mathcal{M}} \in D\bb{\mathcal{M}} \;, \mathcal{M} \in \Omega _T \right\rbrace = \prod _{\mathcal{M} \in \Omega _T} D\bb{\mathcal{M}} \;.
%\end{equation}
\begin{equation}
D\bb{\Omega _T} := \prod _{\mathcal{M} \in \Omega _T} D\bb{\mathcal{M}} \;.
\end{equation}
For any two elements $\bm{u} \;, \bm{v} \in D\bb{\Omega _T}$, and symmetric positive definite permeability tensor $\mathbb{K}$, the weighted $L^2$-inner product is given by
\begin{equation} \label{eq:inner_D}
\bb{\bm{u} \;, \bm{v}}_{\VJ{\mathbb{K}^{-1}}} = \int_ {\Omega} \bm{u}^\transpose \VJ{\mathbb{K}^{-1}} \bm{v}\ \mathrm{d} \Omega = \sum _{i =1}^T \int_ {\mathcal{M}_i} \bm{u}^\transpose \VJ{\mathbb{K}^{-1}} \bm{v} \ \mathrm{d} \mathcal{M} = {{\mathcal{N}}^2 \bb{\bm{u}}}^\transpose {\mathbb{M}}^{(2)}_{\VJ{\mathbb{K}^{-1}}} {\mathcal{N}}^2 \bb{\bm{v}} \;,
\end{equation}
where, $\mathcal{N}^2\bb{\bm{u}}$ and $\mathcal{N}^2\bb{\bm{v}}$ are the column vector of assembled expansion coefficients of all the sub domains and the mass matrix $\mathbb{M}^{(2)}_{\mathbb{K}^{-1}}$ is given as
\begin{equation}
{\mathbb{M}}^{(2)}_{\VJ{\mathbb{K}^{-1}}} = \bs{\begin{array}{ccccc}
\mathbb{M}^{(2)}_{\VJ{\mathbb{K}^{-1}},1} \nl
& \mathbb{M}^{(2)}_{\VJ{\mathbb{K}^{-1}},2} \nl
& & \ddots \nl
& & & \mathbb{M}^{(2)}_{\VJ{\mathbb{K}^{-1}},T}
\end{array}
} \;, \qquad \text{with,} \qquad \mathbb{M}^{(2)}_{\VJ{\mathbb{K}^{-1}},i} = \int _{\mathcal{M}_i} {\Psi^2_i \bb{\bm{x}}}^\transpose {\VJ{\mathbb{K}^{-1}}\bb{\bm{x}} } \Psi ^2_i \bb{\bm{x}} \mathrm{d}\mathcal{M} \;,
\end{equation}
where $\Psi _i ^2 \bb{\bm{x}}$ are the basis functions associated to the finite dimensional space $D\bb{\mathcal{M}_i}$.
\subsubsection{Finite dimensional space $S\bb{\Omega _T} \subset L^2\bb{\Omega _T}$}
For the set of domains $\Omega _T$ we define the finite element space $S\bb{\Omega _T} \subset L^2\bb{\Omega _T}$ as
%\begin{equation}
%S\bb{\Omega _T} := \left\lbrace p |_{\mathcal{M}} \in S \bb{\mathcal{M}} \;, \mathcal{M} \in \Omega _T \right\rbrace = \prod _{\mathcal{M} \in \Omega _T} S \bb{\mathcal{M}} \;.
%\end{equation}
\begin{equation}
S\bb{\Omega _T} := \prod _{\mathcal{M} \in \Omega _T} S \bb{\mathcal{M}} \;.
\end{equation}
For elements $p,q \in S\bb{\Omega _T}$, the $L^2$- inner product is given by
\begin{equation}
\bb{p,q} = \int _{\Omega} p^\transpose q\ \mathrm{d}\Omega = \sum_{i=1}^T  \int _{\mathcal{M}_i} p^\transpose q\ \mathrm{d}\mathcal{M} = {\mathcal{N}^3\bb{p}}^\transpose \mathbb{M}^{(3)} \mathcal{N}^3 \bb{q} \;,
\end{equation}
%\begin{equation}
%\bb{p,q} = \int _{\VJ{\Omega}} p^\transpose q\ \mathrm{d}\Omega = {\mathcal{N}^3\bb{p}}^\transpose \bb{\int _{\VJ{\Omega}} {\mathbf{\Psi}^3\bb{x}}^\transpose \mathbf{\Psi}^3\bb{x}\ \mathrm{d} \Omega} \mathcal{N}^3 \bb{q} = {\mathcal{N}^3\bb{p}}^\transpose \mathbb{M}^{(3)} \mathcal{N}^3 \bb{q} \;,
%\end{equation}
where $\mathcal{N}^3\bb{p}$, $\mathcal{N}^3\bb{q}$ are the column vectors of assembled expansion coefficients of all the sub domains, and the mass matrix $\mathbb{M}^{(3)}$ is the block diagonal mass matrix given by
\begin{equation}
\mathbb{M}^{(3)} = \bs{\begin{array}{ccccc}
\mathbb{M}^{(3)}_1 \nl
& \mathbb{M}^{(3)}_2 \nl
& & \ddots \nl
& & & \mathbb{M}^{(3)}_T
\end{array}
} \;, \qquad \text{with,} \qquad \mathbb{M}^{(3)}_i = \int _{\mathcal{M}_i} {\Psi^3_i \bb{\bm{x}}}^\transpose \Psi ^3_i \bb{\bm{x}} \mathrm{d}\mathcal{M} \;,
\end{equation}
where $\Psi_i ^3\bb{\bm{x}}$ are the basis associated with the finite dimensional space $S\bb{\mathcal{M}_i}$.
%In case of Darcy flow we have the anisotropy tensor $\mathbb{K}$, and therefore we need to evaluate the weighted inner product.
%For domain $\mathcal{M}_i$ the weighted inner product is evaluated as
%\begin{equation}
%{\mathbb{M}^{(2)}_{\mathbb{K}^{-1}}}_i = \int _{\mathcal{M}_i} \Psi_i ^2 \bb{\bm{x}}^\transpose \mathbb{K}^{-1} \bb{\bm{x}}  \Psi_i ^2\bb{\bm{x}} \ \mathrm{d}\mathcal{M} \;.
%\end{equation}
%The weighted inner product term in \eqref{eq:var_hybrid1_finite} is then given by
%\begin{equation}
%\bb{\bm{v} \;, \mathbb{K}^{-1} \bm{u}} = \int _\Omega \bm{v}^\transpose \mathbb{K}^{-1} \bm{u} \ \mathrm{d}\Omega = \sum	_{i=1}^T \int _{\mathcal{M}_i} \bm{v}^\transpose \mathbb{K}^{-1} \bm{u} \ \mathrm{d}\mathcal{M} = \mathcal{N}^2\bb{\bm{v}}^\transpose \bb{\sum _{i=1}^T \int _{\mathcal{M}_i} \Psi_i ^2 \bb{\bm{x}}^\transpose \mathbb{K}^{-1} \bb{\bm{x}}  \Psi_i ^2\bb{\bm{x}} \ \mathrm{d}\mathcal{M}}  \mathcal{N}^2\bb{\bm{u}} = \mathcal{N}^2\bb{\bm{v}}^\transpose \mathbb{M}^{(2)}_{\mathbb{K}^{-1}} \mathcal{N}^2\bb{\bm{u}} \;,
%\end{equation}
%where $\mathcal{N}^2\bb{\bm{v}} \;, \mathcal{N}^2\bb{\bm{u}}$ are the assembled expansion coefficients of all the sub domains, and the mass matrix $\mathbb{M}^{(2)}_{\mathbb{K}^{-1}}$ is given as
%\[\ {\mathbb{M}}_{\mathbb{K}^{-1}}^{(2)} = \bs{\begin{array}{ccccc}
%{\mathbb{M}^{(2)}_{\mathbb{K}^{-1}}}_1 \nl
%& {\mathbb{M}^{(2)}_{\mathbb{K}^{-1}}}_2 \nl
%& & \ddots \nl
%& & & {\mathbb{M}^{(2)}_{\mathbb{K}^{-1}}}_T
%\end{array}} \;. \]
\subsubsection{Finite dimensional space $D_b\bb{\partial \Omega _T} \subset H^{-1/2}\bb{\partial \Omega _T}$}
For the set of boundaries $\partial \Omega_T$ we define the trace space $D_b\bb{\partial \Omega _T} \subset H^{-1/2}\bb{\partial \Omega _T}$ as
\begin{equation}
%D_b\bb{\partial \Omega _T} := \left\lbrace \mu |_{\partial \mathcal{M}} \in D_b\bb{\partial \mathcal{M}} \;, \partial \mathcal{M} \in \partial \Omega _T  \right\rbrace = \prod _{\mathcal{M} \in \Omega _T} D_b\bb{\partial \mathcal{M}} \;.
D_b\bb{\partial \Omega _T} := \prod _{\mathcal{M} \in \Omega _T} D_b\bb{\partial \mathcal{M}} \;.
\end{equation}
\subsubsection{Finite dimensional space $\widetilde{S}\bb{\Omega _T} \subset L^2\bb{\Omega _T}$}
The algebraic dual \spaces of $S\bb{\Omega_T}$ is denoted by $\widetilde{S} \bb{\Omega _T}$  and defined as
%\begin{equation}
%\widetilde{S}\bb{\Omega _T} := \left\lbrace p |_{\mathcal{M}} \in \widetilde{S} \bb{\mathcal{M}} \;, \mathcal{M} \in \Omega _T \right\rbrace = \prod _{\mathcal{M} \in \Omega _T} \widetilde{S} \bb{\mathcal{M}} \;.
%\end{equation}
\begin{equation}
\widetilde{S}\bb{\Omega _T} := \prod _{\mathcal{M} \in \Omega _T} \widetilde{S} \bb{\mathcal{M}} \;.
\end{equation}
%Let $\widetilde{S}\bb{\Omega _T}$ be the dual space of $S\bb{\Omega _T}$.
%If $\mathbf{\widetilde{\Psi}}^0\bb{\bm{x}}$ and $\widetilde{\mathcal{N}}^0\bb{q}$ represent the basis functions and the expansion coefficients, then we have
%\begin{equation}
%\mathbf{\widetilde{\Psi}}^0\bb{\bm{x}} = {\mathbb{M}^{(3)}}^{-1} \mathbf{\Psi}^3 \bb{\bm{x}} \;, \qquad \mbox{and} \;, \qquad \widetilde{\mathcal{N}}^0 \bb{q} = \mathbb{M}^{(3)} \mathcal{N}^3 \bb{q} \;.
%\end{equation}
Let $p,q \in S\bb{\Omega _T} \times \widetilde{S}\bb{\Omega _T}$, then the inner product between the elements is given by
\begin{equation} \label{eq:inner_S_dual}
\bb{p,q} = \sum_{i=1}^T \int _{\mathcal{M}_i} q^\transpose p\ \mathrm{d}\mathcal{M} = {\widetilde{\mathcal{N}}^0 \bb{q}}^\transpose \mathcal{N}^3\bb{p} \;,
\end{equation}
where, $\widetilde{\mathcal{N}}^0 \bb{q} \;, \mathcal{N}^3\bb{p}$ are the column vectors of assembled expansion coefficients.
%\begin{equation} \label{eq:inner_S_dual}
%\VJ{\ba{p,q} :=} \bb{p,q} = {\widetilde{\mathcal{N}}^0 \bb{p}}^\transpose \bb{\int _{\VJ{\Omega}} \mathbf{\widetilde{\Psi}} ^0\bb{\bm{x}}^\transpose \mathbf{\Psi} ^3 \bb{\bm{x}} \mathrm{d}\Omega} \mathcal{N}^3\bb{q} = {\widetilde{\mathcal{N}}^0 \bb{p}}^\transpose \mathcal{N}^3\bb{q} \;.
%\end{equation}
\subsubsection{Finite dimensional space $\widetilde{D}_b\bb{\partial \Omega _e} \subset H^{1/2}\bb{\partial \Omega _e}$}
Let $\gamma \in \partial \Omega _e$ be any boundary face of the sub domains.
Let $D_b\bb{\gamma}$ be the space of restriction of basis of $D\bb{\Omega _T}$ to $\gamma$ and $\widetilde{D}_b\bb{\gamma}$ be its algebraic dual \spaces.

Then the finite dimensional sub space $\widetilde{D}_b\bb{\partial \Omega_e} \subset H^{1/2}\bb{\partial \Omega _e}$ is defined as
\begin{equation}
\widetilde{D}_b\bb{\partial \Omega _e} = \prod _{\gamma \in \partial \Omega _e} \widetilde{D}_b \bb{\gamma} \;.
\end{equation}
%The finite dimensional approximation of the trace space $H^{1/2}\bb{\overline{\partial \Omega _T}}$ requires additional treatment.
%For the interface $\gamma \in \overline{\partial \Omega _T}$ we define the the finite dimensional trace space $D_b\bb{\gamma}$ as
%\[ D_b\bb{\gamma} = \left\lbrace \mu |_{\gamma} \in D_b\bb{\partial \mathcal{M}_i} \quad \mathrm{for} \quad \gamma \subset \partial \mathcal{M}_i \subset \overline{\partial \Omega _T} \right\rbrace \;. \]
%Let $\widetilde{D}_b\bb{\gamma}$ be the corresponding dual space.
%The finite dimensional approximation of $H^{1/2}\bb{\overline{\partial \Omega _T}}$ is then given by
%\[\widetilde{D}_b\bb{\overline{\partial \Omega _T}} = \left\lbrace \widetilde{\lambda} |_{\gamma} \in \widetilde{D}_b\bb{\gamma}, \gamma \in \overline{\partial \Omega _T} \right\rbrace = \prod _{\gamma \in \overline{\partial \Omega _T}} \widetilde{D}_b\bb{\gamma} \;. \] 
For elements $\lambda , \mathrm{tr}\ \bm{u} \in \widetilde{D}_b\bb{\partial \Omega _e} \times D_b\bb{\partial \Omega _T}$, we have the inner product given by
\begin{equation} \label{eq:inner_trace2}
\bb{\lambda,\mathrm{tr}\ \bm{u}} = \sum _{i=1}^T \int _{\partial \mathcal{M}_i} \lambda ^\transpose \mathrm{tr}\ \bm{u}\ \mathrm{d}\Gamma =  \sum _{i=1}^T \widetilde{\mathcal{B}}_i^0 \bb{\lambda}^\transpose \mathcal{B}_i^2 \bb{\mu} \;,
\end{equation} 
where $\widetilde{\mathcal{B}}^0\bb{\lambda}$, $\mathcal{B}^2\bb{\mu}$ are the assembled coefficients over all sub domains.

If $\mu = \bm{u}\cdot \bm{n}|_{\partial \mathcal{M}_i}$ for $i = 1, \hdots , T$, i.e. flux component at sub domain boundaries, then using \eqref{eq:bdof} we have that
\begin{equation} 
\mathcal{B}^2\bb{\mu} = \mathcal{B}^2\bb{\bm{u}\cdot \bm{n}} = \mathbb{N}_2^\transpose \mathcal{N}^2\bb{\bm{u}} \;,
\end{equation}
and we can write \eqref{eq:inner_trace2} as
\begin{equation} \label{eq:inner_trace}
\bb{\lambda , \bm{u} \cdot \bm{n}} = \widetilde{B}^0\bb{\lambda}^\transpose \mathbb{N}_2^\transpose \mathcal{N}^2\bb{\bm{u}} \;,
\end{equation}
where, $\mathcal{N}^2\bb{\bm{u}}$ is the column vector of assembled expansion coefficients and $\mathbb{N}_2^\transpose$ is the assembled trace matrix of all the sub domains.
\subsubsection{The divergence operator for broken Sobolev spaces}
Let $q, \bm{u} \in \widetilde{S}\bb{\Omega _T} \times D\bb{\Omega _T}$ be the elements of the broken finite dimensional spaces.
Then the divergence operation on the vector field $\bm{u}$ is given by
\begin{equation} \label{eq:inner_div}
\bb{q \;, \mathrm{div}\ \bm{u}} = \sum _{i=1}^T \int _{\mathcal{M}_i} q\  \mathrm{div}\ \bm{u}\ \mathrm{d}\mathcal{M} = {\widetilde{\mathcal{N}}^0\bb{q}}^\transpose  \mathbb{E}^{3,2}  \mathcal{N}^2 \bb{\bm{u}}  \qquad \forall q \in \widetilde{S}\bb{\Omega _T} \;,
\end{equation}
%\begin{equation} \label{eq:inner_div}
%\ba{q \;, \mathrm{div}\ \bm{u}} = \bb{q \;, \mathrm{div}\ \bm{u}} = {\widetilde{\mathcal{N}}^0\bb{q}}^\transpose \bb{\int _{\VJ{\Omega}} \mathbf{\widetilde{\Psi}}^0 \bb{\bm{x}}^\transpose \mathbf{\Psi}^3 \bb{\bm{x}}\ \mathrm{d}\Omega } \mathbb{E}^{3,2}  \mathcal{N}^2 \bb{\bm{u}} = {\widetilde{\mathcal{N}}^0\bb{q}}^\transpose  \mathbb{E}^{3,2}  \mathcal{N}^2 \bb{\bm{u}}  \qquad \forall q \in \widetilde{S}\bb{\Omega _T} \;,
%\end{equation}
where, $\widetilde{\mathcal{N}}^0\bb{q}, \mathcal{N}^2\bb{\bm{u}}$ are the column vectors of the assembled expansion coefficients, and $\mathbb{E}^{3,2}$ is the assembled divergence operator.
If $\mathbb{E}^{3,2}_i$ is the discrete representation of divergence operator on the domain $\mathcal{M}_i$, then we have
\begin{equation} \label{eq:incidence_matrix2}
 \mathbb{E}^{3,2} = \bs{ \begin{array}{ccccc}
\mathbb{E}^{3,2}_1 \nl
& \mathbb{E}^{3,2}_2 \nl
& & \ddots \nl
& & & \mathbb{E}^{3,2}_{T}
\end{array}
} \;.
\end{equation}
If the topology of sub domain discretizations is also the same, see \exampleref{ex:topology}, then we have that $\mathbb{E}^{3,2}_1 = \mathbb{E}^{3,2}_2 = \hdots = \mathbb{E}^{3,2}_T$.
In this paper we have only used the case with the same topology for all the sub domains, therefore
\begin{equation}
 \mathbb{E}^{3,2} = \bs{ \begin{array}{ccccc}
\mathbb{E}^{3,2}_i \nl
& \mathbb{E}^{3,2}_i \nl
& & \ddots \nl
& & & \mathbb{E}^{3,2}_i
\end{array}
} \;,
\end{equation}
where $\mathbb{E}^{3,2}_i$ is the topological divergence operator for any of the sub domains.

\subsubsection{The gradient of dual representations}
Let $p$ be a scalar field represented by algebraic dual \spaces $\widetilde{S}\bb{\Omega _T}$, and $\hat{p}$ be the boundary value of the scalar field on the sub domain boundary faces $\gamma \in \partial \Omega_e$.
The gradient operation for dual representations $p \;, \hat{p} \in \widetilde{S}\bb{\Omega _T} \times \widetilde{D}_b\bb{\partial \Omega _e}$ is then defined as, $\widetilde{\mathrm{grad}} : \widetilde{S}\bb{\Omega _T} \times \widetilde{D}_b\bb{\partial\Omega_e} \rightarrow \widetilde{D}\bb{\Omega _T} $, such that
\begin{equation}
\int _\Omega \widetilde{\mathrm{grad}} \bb{p\;,\hat{p} } \bm{u} \ \mathrm{d}\Omega = \sum _{i=1}^T \int _{\mathcal{M}_i} \widetilde{\mathrm{grad}} \bb{p\;,\hat{p}} \bm{u} \ \mathrm{d}\mathcal{M} = \sum _{i=1}^T \bb{- \int _{\mathcal{M}_i} p \bb{\mathrm{div} \ \bm{u}} + \int _{\partial \mathcal{M}_i} \hat{p}\bb{ \bm{u} \cdot \bm{n}} \ \mathrm{d} \Gamma} \quad  \forall \bm{u} \in D\bb{\Omega _T} \;.
\end{equation}
%\begin{equation}
%\widetilde{\mathcal{N}}^1\bb{\widetilde{\mathrm{grad}} \bb{\widetilde{p} \;, \widetilde{\hat{p}}}} = - {\mathbb{E}^{3,2}}^\transpose \widetilde{\mathcal{N}}^0\bb{\widetilde{p}} + \mathbb{N}_2 \widetilde{\mathcal{B}}^0 \bb{\widetilde{\hat{p}}} \;,
%\end{equation}
%where $\widetilde{\mathcal{N}}^0\bb{\widetilde{p}} \;, \widetilde{\mathcal{B}}^0 \bb{\widetilde{\hat{p}}}$ are the column vectors of assembled expansion coefficients, and $\mathbb{E}^{3,2}$, \VJ{see \eqref{eq:incidence_matrix2}}, $\mathbb{N}_2$ are the assembled incidence matrix and the inclusion matrix, respectively.
%
For the pair of elements $p \;, \hat{p} \in \widetilde{S}\bb{\Omega _T} \times \widetilde{D}\bb{\partial \Omega_e}$ the $H^1$-norm is then defined as
\begin{equation} \label{eq:H1_norm}
\| p \|^2_{H^1\bb{\Omega}} = \| p \|^2_{L^2\bb{\Omega}} + \| \widetilde{\mathrm{grad}} \bb{p,\hat{p}} \|^2_{L^2\bb{\Omega}} = \sum _{i=1}^T \left\lbrace \| p \|^2_{L^2\bb{\mathcal{M}_i}} + \| \widetilde{\mathrm{grad}} \bb{p,\hat{p}} \|^2_{L^2\bb{\mathcal{M}_i}} \right\rbrace \;.
\end{equation}
\section{Model anisotropic diffusion problem} \label{sec:model_problem}
In this section we will use the broken finite dimensional spaces defined in \secref{sec:finite_spaces} to derive the algebraic formulation for DD formulation of Darcy problem.

%Let $\Omega \in \mathbb{R}^d$, where $d = \mathrm{dim}\bb{\mathbb{R}}$, be a Lipschitz bounded domain with the domain boundary $\partial \Omega$.
The equations for Darcy problem in the domain $\Omega$ are given by
\begin{equation} \label{eq:Darcy}
\left \lbrace \begin{array}{ll}
\bm{u} + \mathbb{K}\ \mathrm{grad}\ p & = 0 \nl
\mathrm{div}\ \bm{u} & = f
\end{array} \right. \quad \mbox{with} \quad
\left \lbrace \begin{array}{ll}
\bm{u} \cdot \bm{n} = \hat{u} & \mbox{on} \quad \Gamma_N \nl
p = \hat{p} & \mbox{on} \quad \Gamma _D \nl
\Gamma _N \cap \Gamma _D = \emptyset \quad & \mbox{and} \quad \Gamma _N \cup \Gamma _D = \partial \Omega
\end{array} \right. \;,
\end{equation}
where $\bm{u}$ is the velocity, $\mathbb{K}$ is the symmetric positive definite permeability tensor, $p$ is the pressure, $f$ is the given right hand side term, $\bm{n}$ is the outward unit normal vector, $\hat{u}$ is the given velocity boundary condition imposed on Neumann boundary $\Gamma _N$ and $\hat{p}$ is the pressure boundary condition imposed on the Dirichlet boundary $\Gamma _D$.

The Lagrange functional for continuous formulation of Darcy equations is given by
\begin{equation} \label{eq:lag_cont}
\mathcal{L}\bb{\bm{u},p;f,\hat{p},\hat{u}} = \int _\Omega \frac{1}{2} \bm{u}^\transpose \mathbb{K}^{-1} \bm{u} \ \mathrm{d}\Omega - \int _\Omega p \bb{\mathrm{div}\ \bm{u} - f} \mathrm{d}\Omega + \int _{\Gamma _D} \hat{p} \bb{\bm{u}\cdot \bm{n}} \mathrm{d}\Gamma + \int _{\Gamma _N} \bb{\bm{u} \cdot \bm{n} - \hat{u}} \mathrm{d}\Gamma \;.
\end{equation}
The algebraic system for \eqref{eq:lag_cont}, using algebraic dual \representations of \cite{2020jain}, is given by
\begin{equation}
\bs{\begin{array}{cc}
\mathbb{M}_{\mathbb{K}^{-1}}^{(2)} & -{\mathbb{E}^{3,2}}^\transpose \nl
- \mathbb{E}^{3,2} & 0
\end{array}
} \bs{\begin{array}{c}
\mathcal{N}^2\bb{\bm{u}} \nl
\widetilde{\mathcal{N}}^0 \bb{p}
\end{array}
} = \bs{\begin{array}{c}
- {\mathbb{N}_2} \widetilde{\mathcal{B}}^0 \bb{\hat{p}} \nl
- \mathcal{N}^3 \bb{f}
\end{array}
} \;.
\end{equation}
The above system can be solved for unknowns of $p$, and the unknowns of $\bm{u}$ as
\begin{eqnarray}
\widetilde{\mathcal{N}}^0\bb{p} & = & \bb{\mathbb{E}^{3,2} {\mathbb{M}^{(2)}_{\mathbb{K}^{-1}}}^{-1}{\mathbb{E}^{3,2}}^\transpose}^{-1} \bb{\mathcal{N}^3\bb{f} + \mathbb{E}^{3,2} {\mathbb{M}^{(2)}_{\mathbb{K}^{-1}}}^{-1} \mathbb{N}_2 \widetilde{\mathcal{B}}^0\bb{\hat{p}} } \label{eq:cont_p} \nl
\mathcal{N}^2\bb{\bm{u}} & = & {\mathbb{M}_{\mathbb{K}^{-1}}^{(2)}}^{-1} \bb{{\mathbb{E}^{3,2}}^\transpose \widetilde{\mathcal{N}}^0\bb{p} - \mathbb{N}_2 \widetilde{\mathcal{B}}^0\bb{\hat{p}}} \label{eq:cont_u}
\end{eqnarray}
For the DD formulation of \eqref{eq:Darcy} we break the domain $\Omega$ into $T$ sub domains, see \eqref{eq:break_domain}, and we use the Lagrange multipliers $\lambda$ to enforce the required continuity across the sub domains.
The weak formulation for \eqref{eq:Darcy} is then obtained using the Lagrange functional 
%\begin{equation} \label{eq:Lag_broken_domain}
%\mathcal{L} \bb{\bm{u},p, \lambda; f,\hat{p},\hat{u}} = \sum _{i = 1}^T \left\lbrace \int _{\mathcal{M}_i} \frac{1}{2} \bm{u}^\transpose \mathbb{K}^{-1}\bm{u}\ \d \Omega - \int _{\mathcal{M}_i} p \bb{\mathrm{div}\ \uu - f} \d \Omega + \int _{\partial \mathcal{M}_i \setminus (\Gamma _D \cup \Gamma _N)} \lambda \bb{\uu \cdot \nn} \d \Gamma + \int _{\partial \mathcal{M}_i \cap \Gamma _N } \lambda \bb{\bm{u}\cdot \bm{n} - \hat{u}} \d \Gamma + \int _{\partial \mathcal{M}_i \cap \Gamma _D} \hat{p} \bb{\uu \cdot \nn} \d\Gamma \right\rbrace \;,
%\end{equation}
\begin{eqnarray}
\mathcal{L} \bb{\bm{u},p, \lambda; f,\hat{p},\hat{u}} & = & \sum _{i = 1}^T \left\lbrace \int _{\mathcal{M}_i} \frac{1}{2} \bm{u}^\transpose \mathbb{K}^{-1}\bm{u}\ \d \Omega - \int _{\mathcal{M}_i} p \bb{\mathrm{div}\ \uu - f} \d \Omega + \int _{\partial \mathcal{M}_i \setminus (\Gamma _D \cup \Gamma _N)} \lambda \bb{\uu \cdot \nn} \d \Gamma \right. \nonumber \nl
&  & \left. + \int _{\partial \mathcal{M}_i \cap \Gamma _N } \lambda \bb{\bm{u}\cdot \bm{n} - \hat{u}} \d \Gamma + \int _{\partial \mathcal{M}_i \cap \Gamma _D} \hat{p} \bb{\uu \cdot \nn} \d\Gamma \right\rbrace \;, \label{eq:Lag_broken_domain}
\end{eqnarray}
where, on the left hand side, the first term is the kinetic energy term, the second term imposes the constraint on divergence of velocity field $\bm{u}$, the third term imposes continuity of flux across the sub domain faces, the fourth term imposes the Neumann boundary condition and the fifth term imposes the Dirichlet boundary conditions.
%\begin{equation} \label{eq:Lag_broken_domain}
%\mathcal{L} \bb{\bm{u},p, \lambda; f,\hat{p},\hat{u}} = \int _{\VJ{\Omega}} \frac{1}{2} \bm{u}^\transpose \mathbb{K}^{-1}\bm{u}\ \d \Omega - \int _{\VJ{\Omega}} p \bb{\mathrm{div}\ \uu - f} \d \Omega + \int _{\partial \VJ{\Omega}} \lambda \bb{\uu \cdot \nn} \d \Gamma + \int _{\Gamma _N } \lambda \bb{\bm{u}\cdot \bm{n} - \hat{u}} \d \Gamma + \int _{\Gamma _D} \hat{p} \bb{\uu \cdot \nn} \d\Gamma \;,
%\end{equation}
%where, $\lambda$ is the Lagrange multiplier that imposes the constraint on continuity of flux across the sub domains.
\begin{lemma} \label{lem:lambda}
The Lagrange multipliers are pressure boundary values on sub domains.
\begin{proof}
In \eqref{eq:Lag_broken_domain} if we take variations with respect to $\bm{u}$, we get
\[ \sum_{i=1}^T \left\lbrace \int _{\mathcal{M}_i} \bm{v}^\transpose \mathbb{K}^{-1} \bm{u} \ \mathrm{d}\mathcal{M} - \int _{\mathcal{M}_i} \bb{\mathrm{div}\ \bm{v}} p  \ \mathrm{d}\mathcal{M} + \int _{\partial \mathcal{M}_i \setminus {\Gamma _D}} \bb{\bm{v} \cdot \bm{n}} \lambda \ \mathrm{d}\Gamma + \int _{\partial \mathcal{M}_i \cap \Gamma _D} \hat{p} \bb{\bm{v}\cdot \bm{n}} \mathrm{d}\Gamma \right\rbrace = 0 \qquad \forall \bm{v} \in H\bb{\mathrm{div};\Omega _t} \;. \]
%\VJ{where ${\Gamma _D}_i = \partial \mathcal{M}_i \cap \Gamma _D$.}
For any sub domain $\mathcal{M}_i$, we have
\[ \int _{\mathcal{M}_i} \bm{v}^\transpose \mathbb{K}^{-1} \bm{u} \ \mathrm{d}\mathcal{M} - \int _{\mathcal{M}_i} \bb{\mathrm{div}\ \bm{v}} p  \ \mathrm{d}\mathcal{M} + \int _{\partial \mathcal{M}_i \setminus {\Gamma _D}} \bb{\bm{v} \cdot \bm{n}} \lambda \ \mathrm{d}\Gamma + \int _{\partial \mathcal{M}_i \cap \Gamma _D} \bb{\bm{v}\cdot \bm{n}} \hat{p} \mathrm{d}\Gamma = 0 \qquad \forall \bm{v} \in H\bb{\mathrm{div};\mathcal{M}_i} \;. \]
If the solution is sufficiently smooth, then using integration by parts on the second term, we get
\[ \int _{\mathcal{M}_i} \bm{v}^\transpose \mathbb{K}^{-1} \bm{u} \ \mathrm{d}\mathcal{M} - \bb{- \int_{\mathcal{M}_i} \bm{v}^\transpose \bb{\mathrm{grad}\ p} \mathrm{d}\mathcal{M} + \int _{\partial \mathcal{M}_i} \bb{\bm{v}\cdot \bm{n}} p \ \mathrm{d}\Gamma } + \int _{\partial \mathcal{M}_i \setminus {\Gamma _D}} \bb{\bm{v} \cdot \bm{n}} \lambda \ \mathrm{d}\Gamma + \int _{\partial \mathcal{M}_i \cap \Gamma _D} \bb{\bm{v}\cdot \bm{n}} \hat{p} \ \mathrm{d}\Gamma = 0 \quad \forall \bm{v} \in H\bb{\mathrm{div};\mathcal{M}_i} \;. \]
Now, combining the first and the second term, and the third, fourth and fifth term we get
\[ \int _{\mathcal{M}_i} \bm{v}^\transpose \bb{\mathbb{K}^{-1} \bm{u} + \mathrm{grad} \ p} \ \mathrm{d}\mathcal{M} + \int _{\partial \mathcal{M}_i \setminus {\Gamma _D}} \bb{\bm{v} \cdot \bm{n}} \bb{\lambda - p} \mathrm{d}\Gamma + \int _{\partial \mathcal{M}_i \cap \Gamma _D} \bb{\bm{v}\cdot\bm{n}} \bb{p - \hat{p}} \mathrm{d}\Gamma = 0 \qquad \forall \bm{v} \in H\bb{\mathrm{div};\mathcal{M}_i} \;. \] 
This implies that $\mathbb{K}^{-1} \bm{u} + \mathrm{grad}\ p =0$ in the $L^2$-sense, i.e. almost everywhere, and $\lambda -p =0$ in the $H^{1/2}$-sense, i.e. between the sub domains and $p = \hat{p}$ along $\Gamma _D$.
As this should hold for all $\bm{v} \in H\bb{\mathrm{div};\mathcal{M}_i}$, we have that $\lambda = p$, i.e. the Lagrange multipliers are the pressure boundary values on the sub domain boundaries, $\partial \mathcal{M}_i$.
\end{proof}
\end{lemma}
The optimality conditions at continuous level	for the Lagrange functional \eqref{eq:Lag_broken_domain} are given by: For given $f \in L^2\bb{\Omega}$, $\hat{p} \in H^{1/2} \bb{\Gamma _D}$, $\hat{u} \in H^{-1/2}\bb{\Gamma _N}$, find $\bm{u} \in H\bb{\mathrm{div};\mathcal{M}_i}$, $p \in L^2\bb{\mathcal{M}_i}$, $\lambda \in H^{1/2}\bb{\partial \Omega _e}$, such that 
\begin{equation} \label{eq:var_hybrid1_cont}
\sum _{i=1}^T
\left\lbrace
\begin{array}{clllll}
\bb{\bm{v}, \mathbb{K}^{-1} \uu} & - \bb{\mathrm{div}\ \bm{v}, p} & + \bb{\bm{v} \cdot \nn, \lambda} & = -\bb{\bm{v}\cdot \nn, \hat{p}} & \forall \bm{v} \in H\bb{\mathrm{div}; \mathcal{M}_i} \nl
-\bb{q, \mathrm{div}\ \uu} & & & = -\bb{q, f} & \forall q \in L^2\bb{\mathcal{M}_i} \nl
\bb{\mu , \uu \cdot \nn} & & & = \bb{\mu ,\hat{u}} & \forall \mu \in H^{1/2} \bb{\partial \Omega _e}
\end{array}
\right. \;.
\end{equation}
The finite dimensional problem is then given by: For given $f \in L^2\bb{\Omega}$, $\hat{p} \in H^{1/2} \bb{\Gamma _D}$, $\hat{u} \in H^{-1/2}\bb{\Gamma _N}$, find $\bm{u} \in D\bb{\mathcal{M}_i}$, $ {p} \in \widetilde{S}\bb{\mathcal{M}_i}$, ${\lambda} \in \widetilde{D}_b\bb{\partial \Omega _e}$, such that 
\begin{equation} \label{eq:var_hybrid1_finite}
\sum _{i=1}^T
\left\lbrace
\begin{array}{clllll}
\bb{\bm{v}, \mathbb{K}^{-1} \uu} & - \bb{\mathrm{div}\ \bm{v}, {p}} & + \bb{\bm{v} \cdot \nn, {\lambda}} & = -\bb{\bm{v}\cdot \nn, {\hat{p}}} & \forall \bm{v} \in D\bb{\mathcal{M}_i} \nl
-\bb{{q}, \mathrm{div}\ \uu} & & & = -\bb{{q}, f} & \forall q \in \widetilde{S}\bb{\mathcal{M}_i} \nl
\bb{{\mu} , \uu \cdot \nn} & & & = \bb{{\mu} ,\hat{u}} & \forall {\mu} \in \widetilde{D}_b \bb{\partial \Omega _e}
\end{array}
\right. \;.
\end{equation}
In \eqref{eq:var_hybrid1_finite} we see that all terms, except one - the weighted inner product term, are the inner product between the primal and the dual representations.
Consequently, these terms do not require evaluation of (dense) mass matrices.
The only matrices associated with these terms are the sparse, metric-free divergence operator, or the inclusion matrix.
See also \eqref{eq:alg_hybrid1}.
%\VJ{
%\subsection{Uniqueness and stability of discrete formulation}
%!!! TBD, to add or not to add, uniqueness and stability of discrete formulation ??? !!! }
\subsection{The algebraic formulation}
The inner product for the left hand side terms of \eqref{eq:var_hybrid1_finite} are evaluated as
\begin{eqnarray}
\bb{\bm{v}, \mathbb{K}^{-1} \uu} & \stackrel{\eqref{eq:inner_D}}{=} & {\mathcal{N}^2\bb{\bm{v}}}^\transpose \mathbb{M}^{(2)}_{\mathbb{K}^{-1}} \mathcal{N}^2\bb{\bm{u}} \;, \nonumber \nl
\bb{{q} \;, \mathrm{div}\ \bm{u}} & \stackrel{\eqref{eq:inner_div}}{=} & {\widetilde{\mathcal{N}}^0 \bb{{q}}}^\transpose \mathbb{E}^{3,2} \mathcal{N}^2\bb{\bm{u}} \;, \nonumber \nl
\bb{{\mu} \;, \bm{u} \cdot \bm{n}} & \stackrel{\eqref{eq:inner_trace}}{=} & {\widetilde{\mathcal{B}}^0 \bb{{\mu}}}^\transpose \mathbb{N}_2^\transpose \mathcal{N}^2\bb{\bm{u}} \nonumber \;,
\end{eqnarray}
and the inner product for the right hand side terms of \eqref{eq:var_hybrid1_finite} are evaluated as
\begin{eqnarray}
\bb{\bm{v} \cdot \bm{n} \;, {\hat{p}}} & \stackrel{\eqref{eq:inner_trace}}{=} & {\mathcal{N}^2\bb{\bm{v}}}^\transpose \mathbb{N}_2 \widetilde{\mathcal{B}}^0 \bb{{\hat{p}}} \VJ{= {\mathcal{N}^2\bb{\bm{v}}}^\transpose \mathbb{N}_2 \int _{\Gamma _D}\widetilde{\Psi}^0_b \bb{\bm{x}} \hat{p}\bb{\bm{x}} \mathrm{d}\Gamma  } \;, \nonumber \nl
\bb{{q} \;, f} & \stackrel{\eqref{eq:inner_S_dual}}{=} & {\widetilde{\mathcal{N}}^0\bb{{q}}}^\transpose \mathcal{N}^3\bb{f} \VJ{= {\widetilde{\mathcal{N}}^0\bb{{q}}}^\transpose \int _\Omega \widetilde{\Psi}^0 \bb{\bm{x}} f\bb{x} \mathrm{d}\Omega }  \;, \nonumber \nl
\bb{{\mu} \;, \hat{u}} & \stackrel{\eqref{eq:inner_trace2}}{=} & {\widetilde{\mathcal{B}}^0 \bb{{\mu}}}^\transpose \mathcal{B}^2 \bb{\hat{u}} \VJ{= {\widetilde{\mathcal{B}}^0 \bb{{\mu}}}^\transpose \int _{\Gamma _N}\Psi_b^2 \bb{\bm{x}} \hat{u} \bb{\bm{x}} \mathrm{d} \Gamma } \nonumber \;.
\end{eqnarray}
Using above relations we can write the algebraic formulation for \eqref{eq:var_hybrid1_finite} as
\begin{equation} \label{eq:alg_hybrid1}
\bs{
\begin{array}{lcc}
\mathbb{M}^{(2)}_{\mathbb{K}^{-1}}& -{\mathbb{E}^{3,2}}^\transpose & \mathbb{N}_2 \nl
-\mathbb{E}^{3,2} & 0 & 0 \nl
\mathbb{N}_2^\transpose & 0 & 0
\end{array}
}\bs{\begin{array}{c}
\mathcal{N}^2\bb{\bm{u}} \nl
\widetilde{\mathcal{N}}^0 \bb{{p}} \nl
\widetilde{\mathcal{B}}^0 \bb{{\lambda}}
\end{array}
} = \bs{\begin{array}{c}
- \mathbb{N}_2 \widetilde{\mathcal{B}}^0\bb{{\hat{p}}} \nl
-\mathcal{N}^3\bb{f} \nl
\mathcal{B}^2 \bb{\hat{u}}
\end{array}
} \;.
\end{equation}
In \eqref{eq:alg_hybrid1}, we see that the matrices $\mathbb{E}^{3,2}$, $\mathbb{N}_2$, are sparse and metric-free.
By metric-free, we mean independent of the size of the elements, the shape of the elements (orthogonal or highly curved) or the order of the approximation.
All the metric dependence is contained in the basis functions and therefore in the mass matrix $\mathbb{M}_{\mathbb{K}^{-1}}^{(2)}$.
%The only dense and metric dependent matrix is $\mathbb{M}^{(2)}_{\mathbb{K}^{-1}}$.
Using static condensation \eqref{eq:alg_hybrid1} can be solved efficiently for the trace variables $\widetilde{\mathcal{B}}^0\bb{{\lambda}}$ only.
We first solve for a global system of $\widetilde{\mathcal{B}}^0\bb{{\lambda}}$ by
%\begin{equation}
%\mathbb{N} \bb{\mathbb{M}^{-1} - \mathbb{M}^{-1} \mathbb{E}^\transpose \bb{\mathbb{E} \mathbb{M}^{-1} \mathbb{E}^\transpose}^{-1} \mathbb{E}\mathbb{M}^{-1}} \mathbb{N}^\transpose \widetilde{\mathcal{B}}^0 \bb{\lambda} = - \mathbb{N}_{\hat{u}} \mathcal{B}^2\bb{\hat{u}} + \mathbb{N} \bb{\mathbb{M}^{-1} - \mathbb{M}^{-1} \mathbb{E}^\transpose \bb{\mathbb{E} \mathbb{M}^{-1} \mathbb{E}^\transpose}^{-1} \mathbb{E}\mathbb{M}^{-1}} \mathbb{N}^\transpose
%\end{equation}
%\[\bb{\conn B_1 \conn^\transpose} \widetilde{\mathcal{N}}_b^0 \bb{\lambda} = - \mathbb{N} _{\hat{u}}\mathcal{B}^2 \bb{\hat{u}} + \conn \bb{\mathbb{M}^{-1} \mathbb{E}^\transpose \bb{\mathbb{E}\mathbb{M}^{-1} \mathbb{E}^\transpose}^{-1} \mathcal{N}^3 \bb{f} - B_1 \mathbb{N}_{\hat{p}}\widetilde{\mathcal{B}}^0 \bb{\hat{p}} } \]
\begin{equation} \label{eq:lambda}
\bb{\conn \mathbb{A} \mathbb{N}_2} \widetilde{\mathcal{B}}^0\bb{{\lambda}} = - \mathcal{B}^2 \bb{\hat{u}} + \conn \bb{{\mathbb{M}_{\mathbb{K}^{-1}}^{(2)}}^{-1} {\mathbb{E}^{3,2}}^\transpose \bb{{\mathbb{E}^{3,2}}{\mathbb{M}_{\mathbb{K}^{-1}}^{(2)}}^{-1} {\mathbb{E}^{3,2}}^\transpose}^{-1} \mathcal{N}^3 \bb{f} - \mathbb{A} \mathbb{N}_2 \widetilde{\mathcal{B}}^0 \bb{{\hat{p}}} }
\end{equation}
%\[\bb{\conn A \conn^\transpose} \widetilde{\mathcal{B}}^0\bb{\lambda} = - \mathbb{N} _{\hat{u}}\mathcal{B}^2 \bb{\hat{u}} + \conn \bb{{\mathbb{M}^{(2)}}^{-1} {\mathbb{E}^{3,2}}^\transpose \bb{{\mathbb{E}^{3,2}}{\mathbb{M}^{(2)}}^{-1} {\mathbb{E}^{3,2}}^\transpose}^{-1} \mathcal{N}^3 \bb{f} - A \mathbb{N}_{\hat{p}}\widetilde{\mathcal{B}}^0 \bb{\hat{p}} } \]
where,
\[\mathbb{A} = {\mathbb{M}_{\mathbb{K}^{-1}}^{(2)}}^{-1} - {\mathbb{M}_{\mathbb{K}^{-1}}^{(2)}}^{-1} {\mathbb{E}^{3,2}}^\transpose \bb{{\mathbb{E}^{3,2}}{\mathbb{M}_{\mathbb{K}^{-1}}^{(2)}}^{-1} {\mathbb{E}^{3,2}}^\transpose}^{-1} {\mathbb{E}^{3,2}} {\mathbb{M}_{\mathbb{K}^{-1}}^{(2)}}^{-1} \;. \]
The matrix $\mathbb{A}$ is block diagonal and can be constructed efficiently by evaluating for each sub domain separately as
\begin{equation} \label{eq:A_i}
\mathbb{A}_i = {\mathbb{M}_{\mathbb{K}^{-1},i}^{(2)}}^{-1} - {\mathbb{M}_{\mathbb{K}^{-1},i}^{(2)}}^{-1} {\mathbb{E}_i^{3,2}}^\transpose \bb{{\mathbb{E}_i^{3,2}}{\mathbb{M}_{\mathbb{K}^{-1},i}^{(2)}}^{-1} {\mathbb{E}_i^{3,2}}^\transpose}^{-1} {\mathbb{E}_i^{3,2}} {\mathbb{M}_{\mathbb{K}^{-1},i}^{(2)}}^{-1} \;.
\end{equation}
The inverse matrices in \eqref{eq:A_i} are symmetric, positive definite and are evaluated using Cholesky decomposition.
The left hand side of \eqref{eq:lambda} is also symmetric, positive definite and the $\lambda$ system is solved using Cholesky decomposition.
%\begin{eqnarray}
%A_1 & = & \mathbb{M}^{-1} \mathbb{E}^\transpose \bb{\mathbb{E}\mathbb{M}^{-1} \mathbb{E}^\transpose}^{-1} \mathcal{N}^3 \bb{f} \nonumber \nl
%B_1 & = & \mathbb{M}^{-1} - \mathbb{M}^{-1} \mathbb{E}^\transpose \bb{\mathbb{E}\mathbb{M}^{-1} \mathbb{E}^\transpose}^{-1} \mathbb{E} \mathbb{M}^{-1} \;, \nonumber
%\end{eqnarray}
The local degrees of freedom in sub domains $\mathcal{M}_i$, $i = 1, \hdots , T$ are then evaluated as
\begin{eqnarray}
\widetilde{\mathcal{N}}^0_i \bb{{p}} & = & \bb{\mathbb{E}^{3,2}_i {\mathbb{M}_{\mathbb{K}^{-1},i}^{(2)}}^{-1}{\mathbb{E}^{3,2}_i}^\transpose}^{-1} \bb{\mathcal{N}_i^3\bb{f} + \mathbb{E}^{3,2}_i {\mathbb{M}_{\mathbb{K}^{-1},i}^{(2)}}^{-1} {\mathbb{N}_2}_i \bb{{ \widetilde{\mathcal{B}}_i^0\bb{{\hat{p}}} + \widetilde{\mathcal{B}}_i^0 \bb{{\lambda}} }} } \label{eq:velocity} \nl
\mathcal{N}^2_i \bb{\bm{u}} & = & {\mathbb{M}_{\mathbb{K}^{-1},i}^{(2)}}^{-1} \bb{{\mathbb{E}^{3,2}_i}^\transpose \widetilde{\mathcal{N}}_i^0 \bb{{p}} - {\mathbb{N}_2}_i  \bb{\widetilde{\mathcal{B}}^0 \bb{\hat{p}}_i + \widetilde{\mathcal{B}}_i^0\bb{{\lambda}}} } \label{eq:pressure} \;.
\end{eqnarray}
The inverse terms in \eqref{eq:velocity} and \eqref{eq:pressure} are already evaluated in \eqref{eq:A_i}, and therefore evaluation of expansion coefficients of velocity and pressure field is simply a matrix multiplication step. 
%where the subscript $i$ denotes that the matrix is associated to the sub domain $\mathcal{M}_i$.
\section{Test cases} \label{sec:test_cases}
In this section we present the computational results for two test cases using the DD formulation.
All the simulations are executed on MATLAB release 2020b on a Macintosh machine with 2.6 GHz Intel Core i7 processor using a single processor.
\subsection{Test case I: Manufactured solution} \label{sec:test_cases_manu}
In this section we will solve a test problem from \cite{2012Wheeler} and compare the results of DD formulation with the continuous formulation for hexahedral elements of order $N=1,2,3$ with varying mesh refinements.
We will show that i) the results from both the formulation are same up to machine precision, ii) that the DD formulation has optimal convergence rates, iii) the speed-up in simulation run times using DD formulation, \VJ{and iv) comparison of condition number of the \eqref{eq:cont_p} and \eqref{eq:lambda}}.
%present the convergence results and the efficiency gained in terms of solution time using the domain decomposition scheme of \eqref{eq:lambda} - \eqref{eq:pressure}, 
%We will also compare the time required to solve this problem using a continuous formulation \cite{2020jain} and DD method described in \secref{sec:model_problem}.
\begin{figure}
\centering
\includegraphics[scale=0.4]{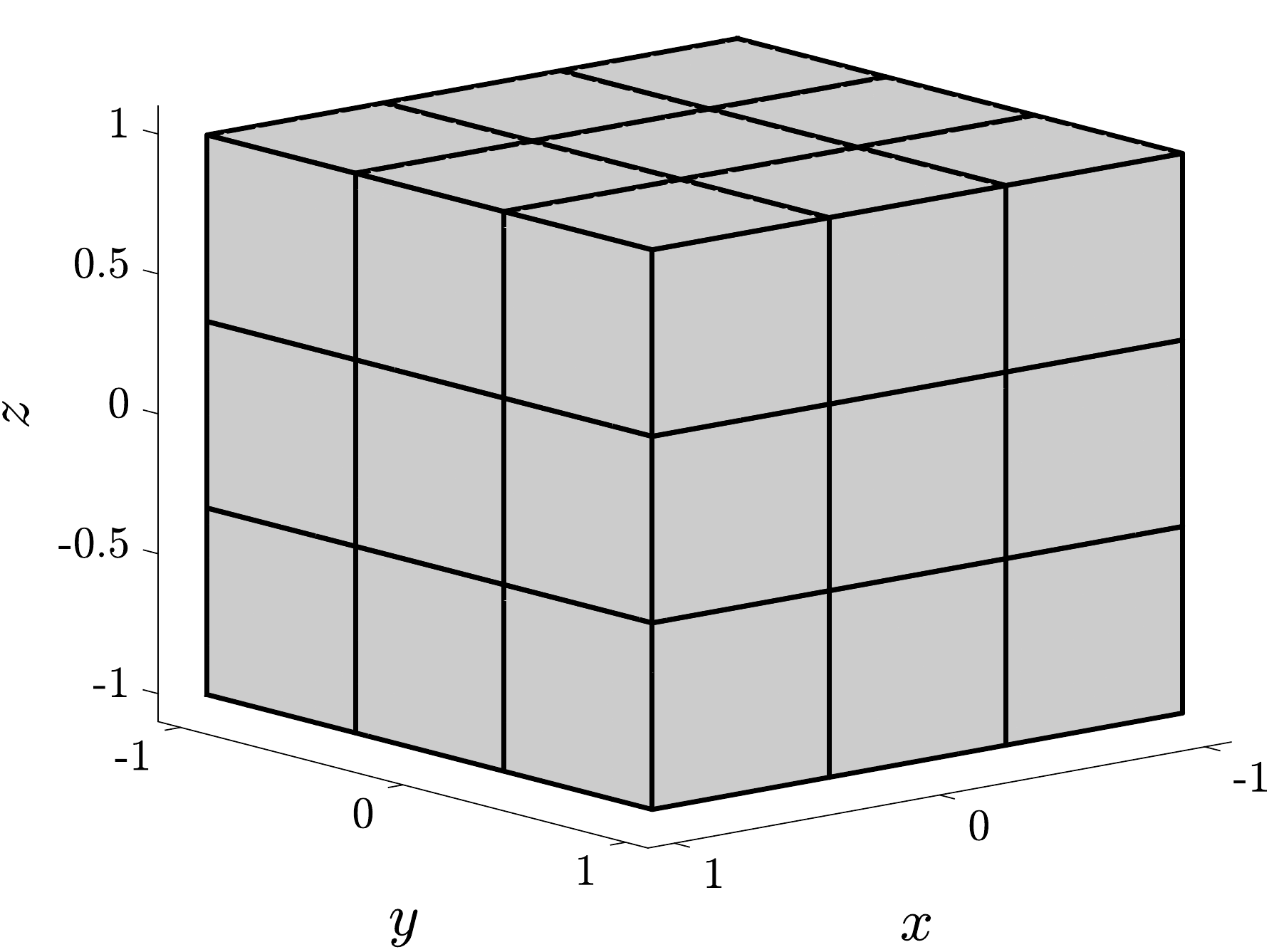} \qquad
\includegraphics[scale=0.4]{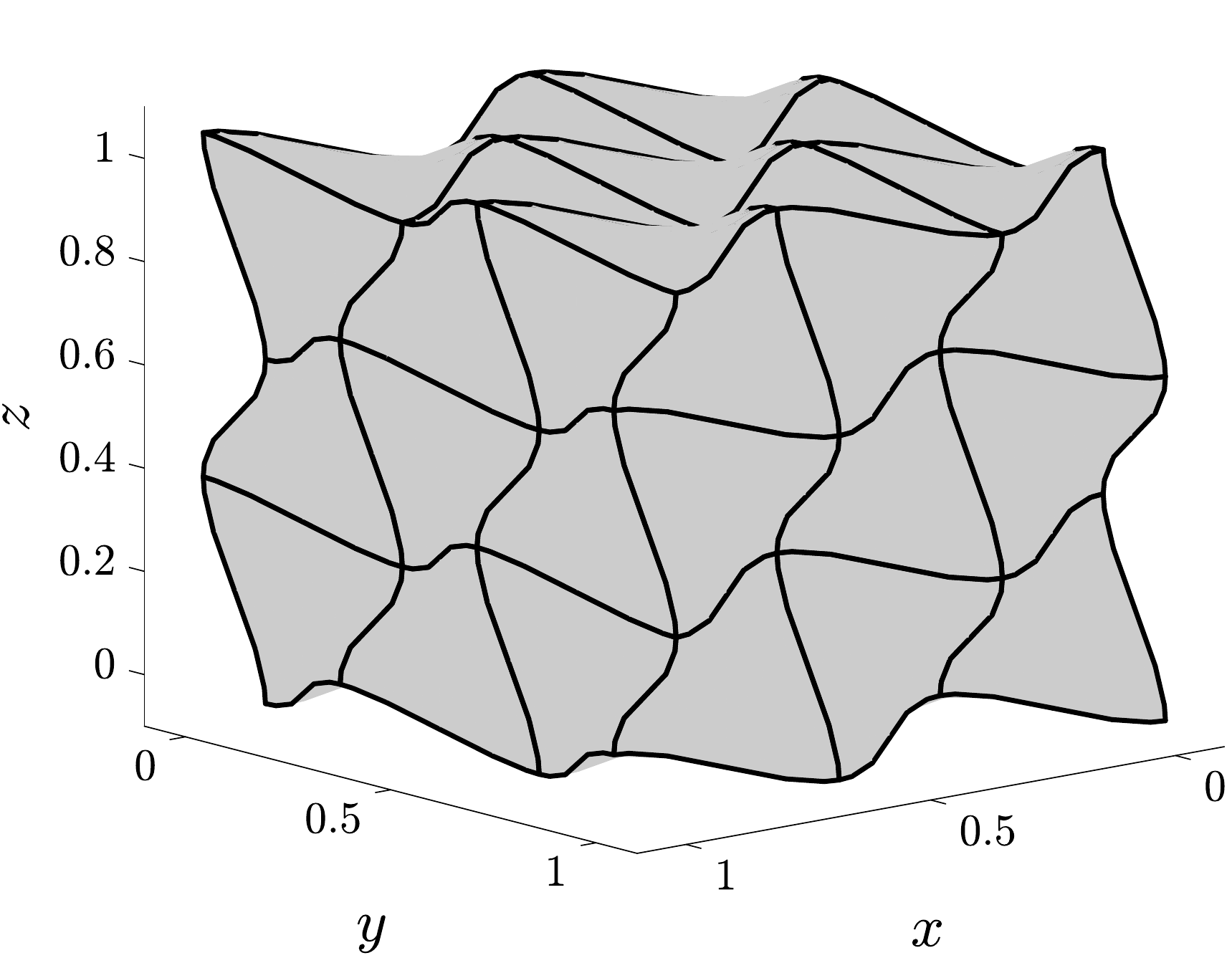}
\caption{Test domain with total number of elements = $3 \times 3 \times 3$. Left: reference domain, $\widehat{\Omega}$, Right: domain of the test case, $\Omega$.}
\label{fig:domain_3D}
\end{figure}

The domain for the problem is obtained by mapping the reference domain $\bb{\xi \;, \eta \;, \zeta} \in \widehat{\Omega} = [-1,1]^3$ by
\begin{equation} \label{eq:mesh}
\left\lbrace
\begin{array}{l}
x = \hat{x} + 0.03 \cos \bb{3 \pi \hat{x}} \cos \bb{3 \pi \hat{y}} \cos \bb{3 \pi \hat{z}}  \\[1.5ex]
y = \hat{y} - 0.04 \cos \bb{3 \pi \hat{x}} \cos \bb{3 \pi \hat{y}} \cos \bb{3 \pi \hat{z}}  \\[1.5ex]
z = \hat{z} + 0.05 \cos \bb{3 \pi \hat{x}} \cos \bb{3 \pi \hat{y}} \cos \bb{3 \pi \hat{z}}
\end{array} \right. \quad \mbox{where,} \quad 
\left\lbrace
\begin{array}{l}
\hat{x} = 0.5 \bb{1 + \xi}  \\[1.5ex]
\hat{y} = 0.5 \bb{1 + \eta}  \\[1.5ex]
\hat{z} = 0.5 \bb{1 + \zeta}
\end{array} \right.  \;.
\end{equation}
In \figref{fig:domain_3D}, in the left plot we show the reference domain $\widehat{\Omega}$ and on the right plot we show the domain of the problem which is obtained using \eqref{eq:mesh}.

The permeability tensor, $\mathbb{K}$, and the exact solution, $p_{ex}$, are given by
\begin{equation*}
\mathbb{K} \bb{x,y,z} = \bs{\begin{array}{ccc}
x^2 + y^2 +1 & 0 & 0 \nl
0 & z^2 +1 & \sin \bb{xy} \nl
0 & \sin \bb{xy} & x^2 y^2 + 1
\end{array}} \;, \qquad \mbox{and} \qquad p_{ex}\bb{x,y,x} = x + y + z - 1.5 \;,
\end{equation*}
and the right hand side term is given by
\begin{equation*}
f_{ex} = -\mathrm{div}\ \mathbb{K} \mathrm{grad} p_{ex} \;.
\end{equation*}
As in \cite{2012Wheeler}, we impose the Dirchlet boundary conditions $\hat{p}$ at $\hat{x}=0,1$ faces and Neumann boundary conditions $\hat{u}$ at $\hat{y} = 0,1, \hat{z}=0,1$ faces, where
\begin{equation*}
\hat{p} = p_{ex}  | _{\Gamma _D} \qquad \mbox{and} \qquad \hat{u} = \bb{ - \mathbb{K} \mathrm{grad} p_{ex} \cdot \bm{n}} |_{\Gamma _N} \;.
\end{equation*}

For DD formulation of this test case we decompose the domain $\Omega$ with equal number of sub domains, $K1$, in each direction.
%The length of elements, $h$, is equal for the unmapped domain and is given by $h = 1/K$
For $N=1,2$ cases each of the sub domains is discretized into $2 \times 2 \times 2$ elements.
For $N=3$ case each sub domain consists of a single element only.
Therefore, total number of elements in one direction, $K$, are $2K1$ for $N=1,2$ case, and $K1$ for $N=3$ case.
We define $\widehat{h} = 2/K$ as the size of a non-deformed element of the reference domain $\widehat{\Omega} = [-1,1]^3$.

\begin{figure}
\includegraphics[width=0.33\linewidth]{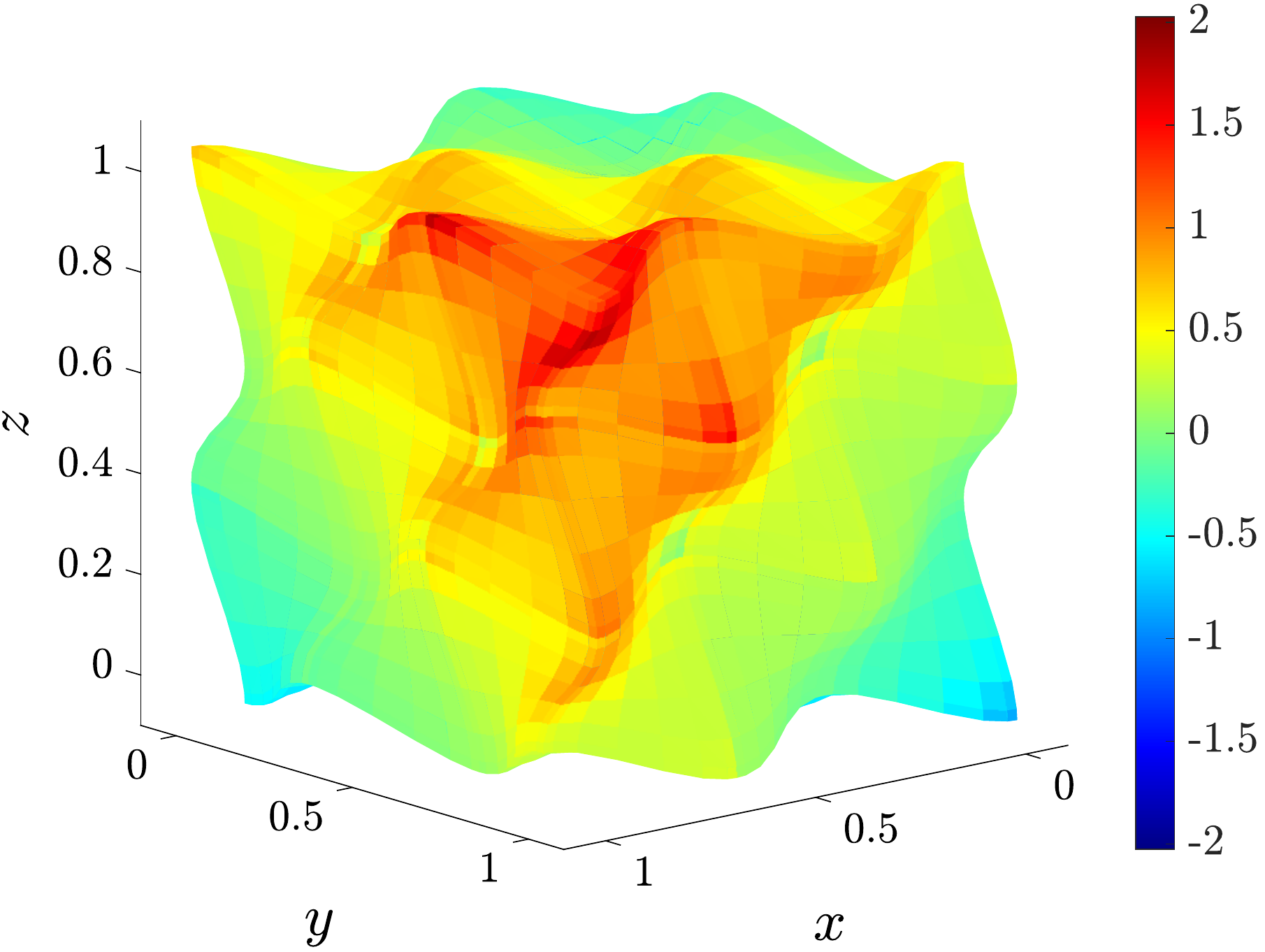}
\includegraphics[width=0.33\linewidth]{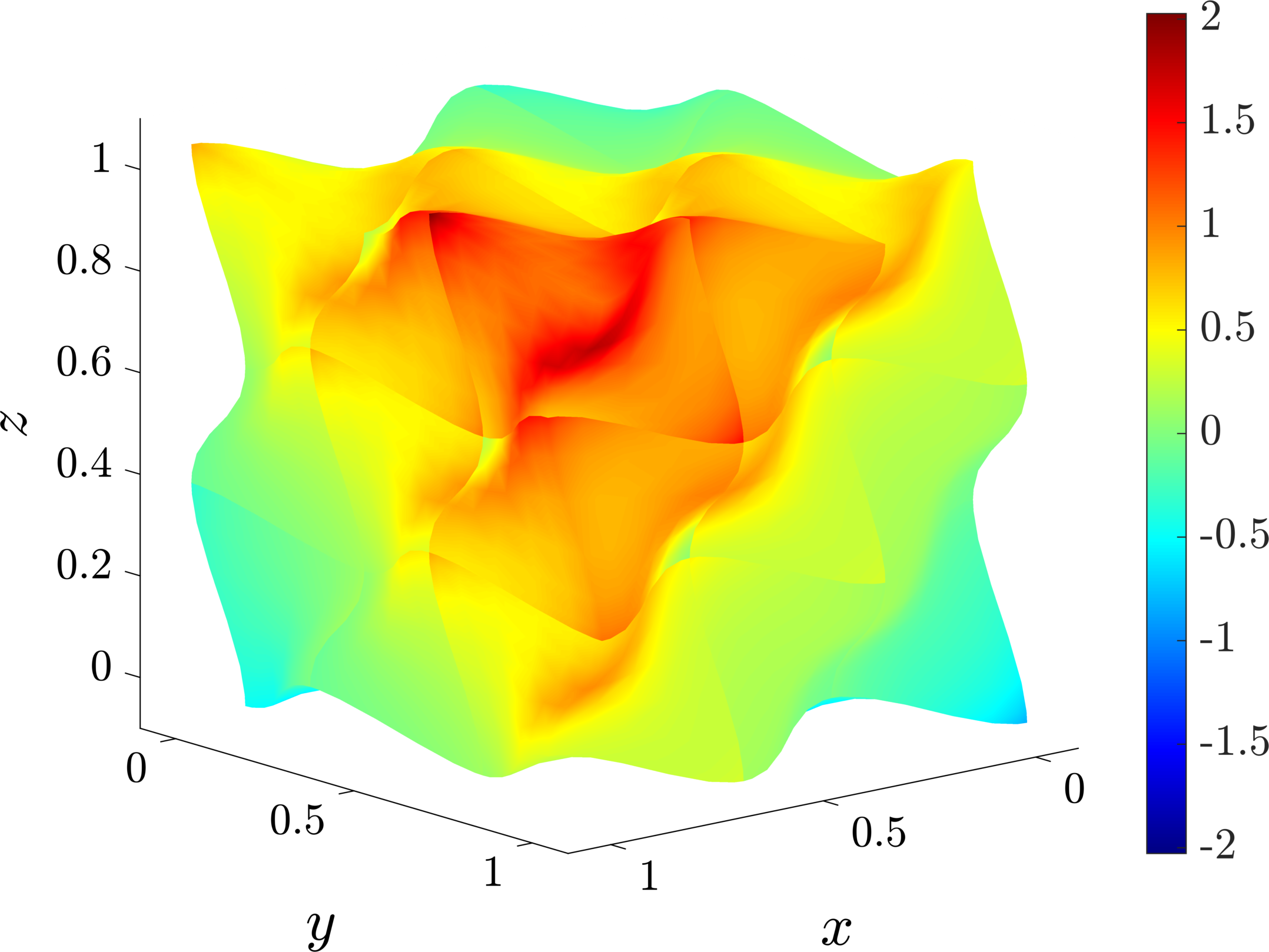}
\includegraphics[width=0.33\linewidth]{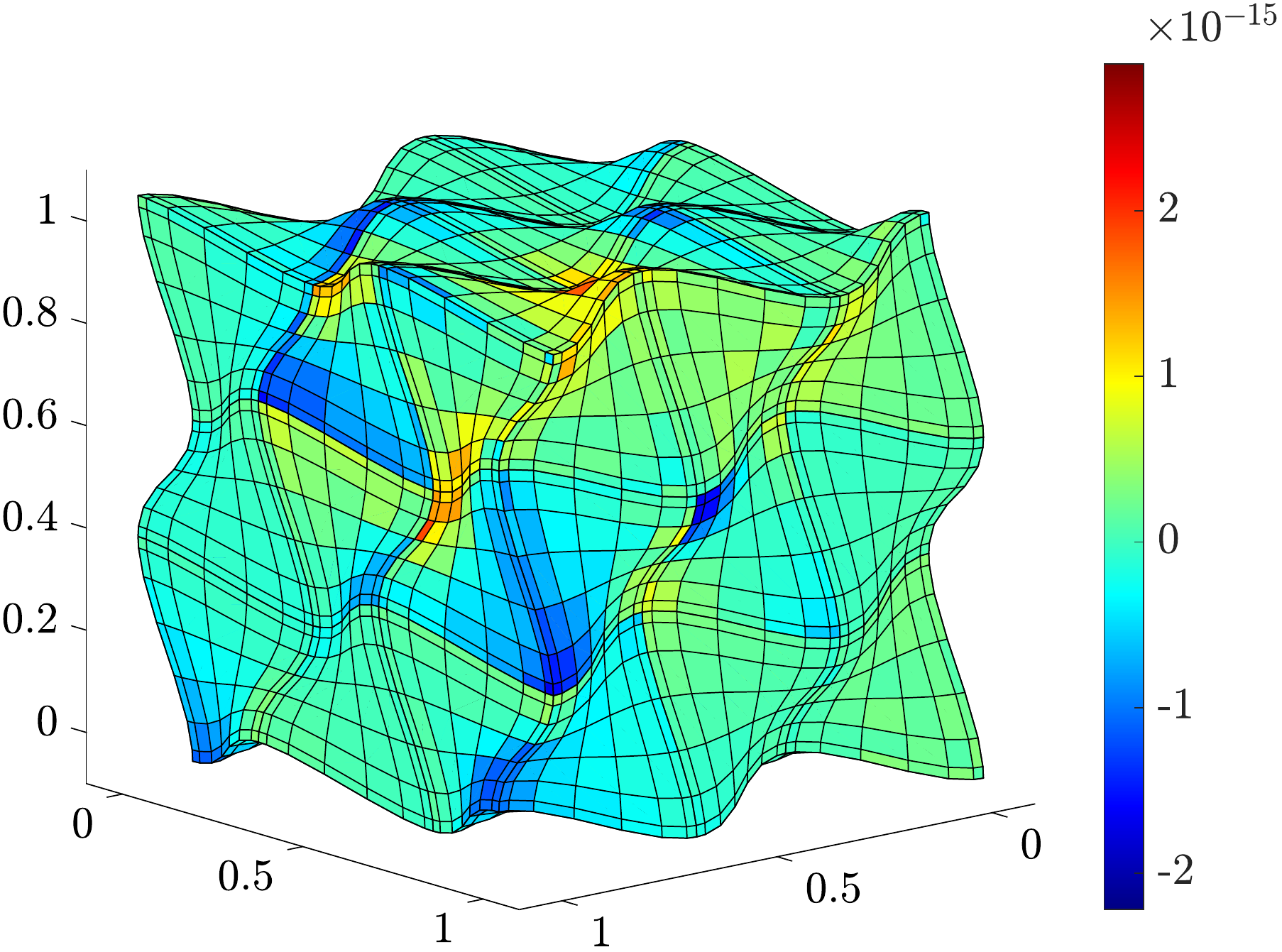} \nl
\includegraphics[width=0.33\linewidth]{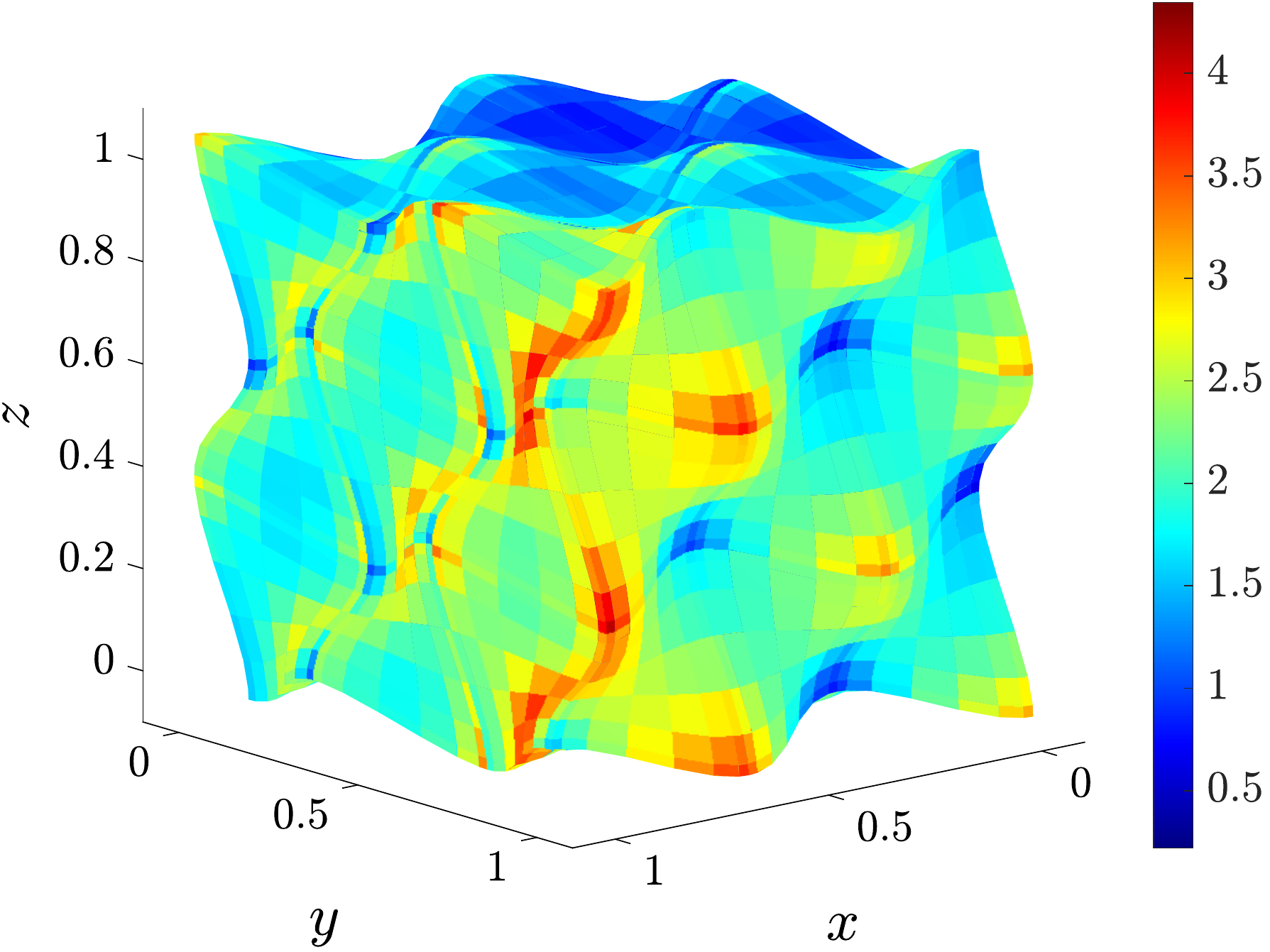}
\includegraphics[width=0.33\linewidth]{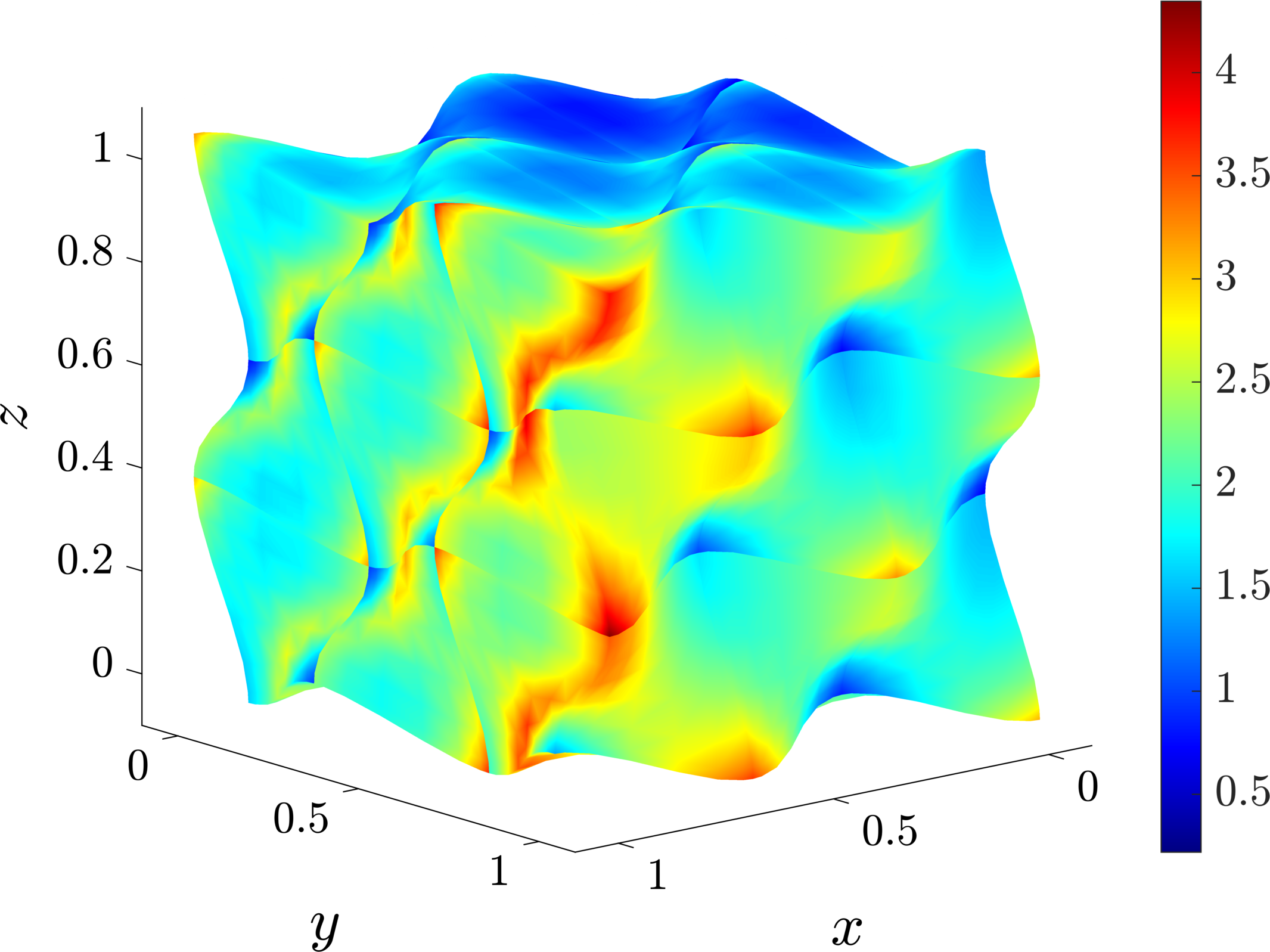}
\includegraphics[width=0.33\linewidth]{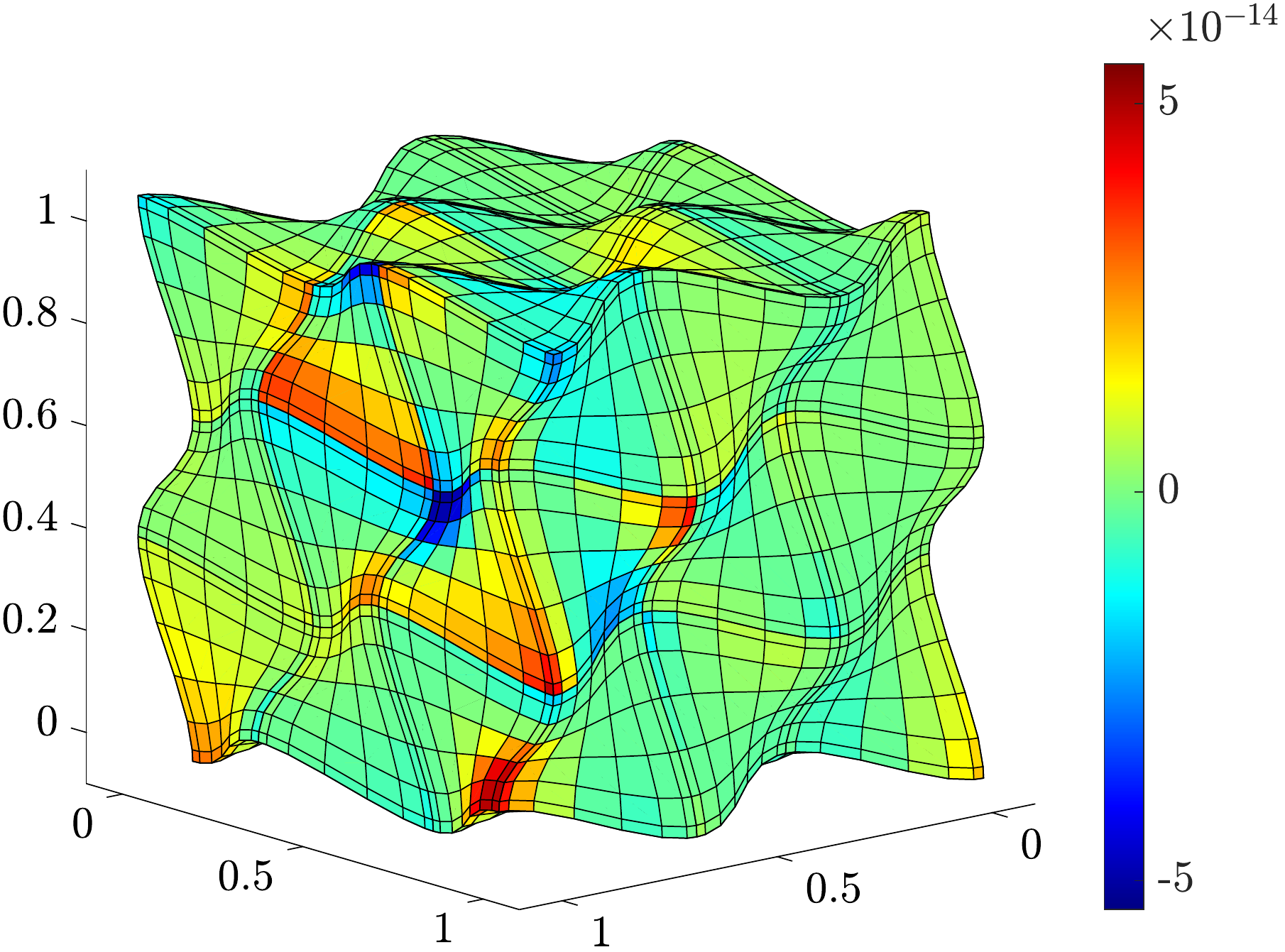} \nl
\includegraphics[width=0.33\linewidth]{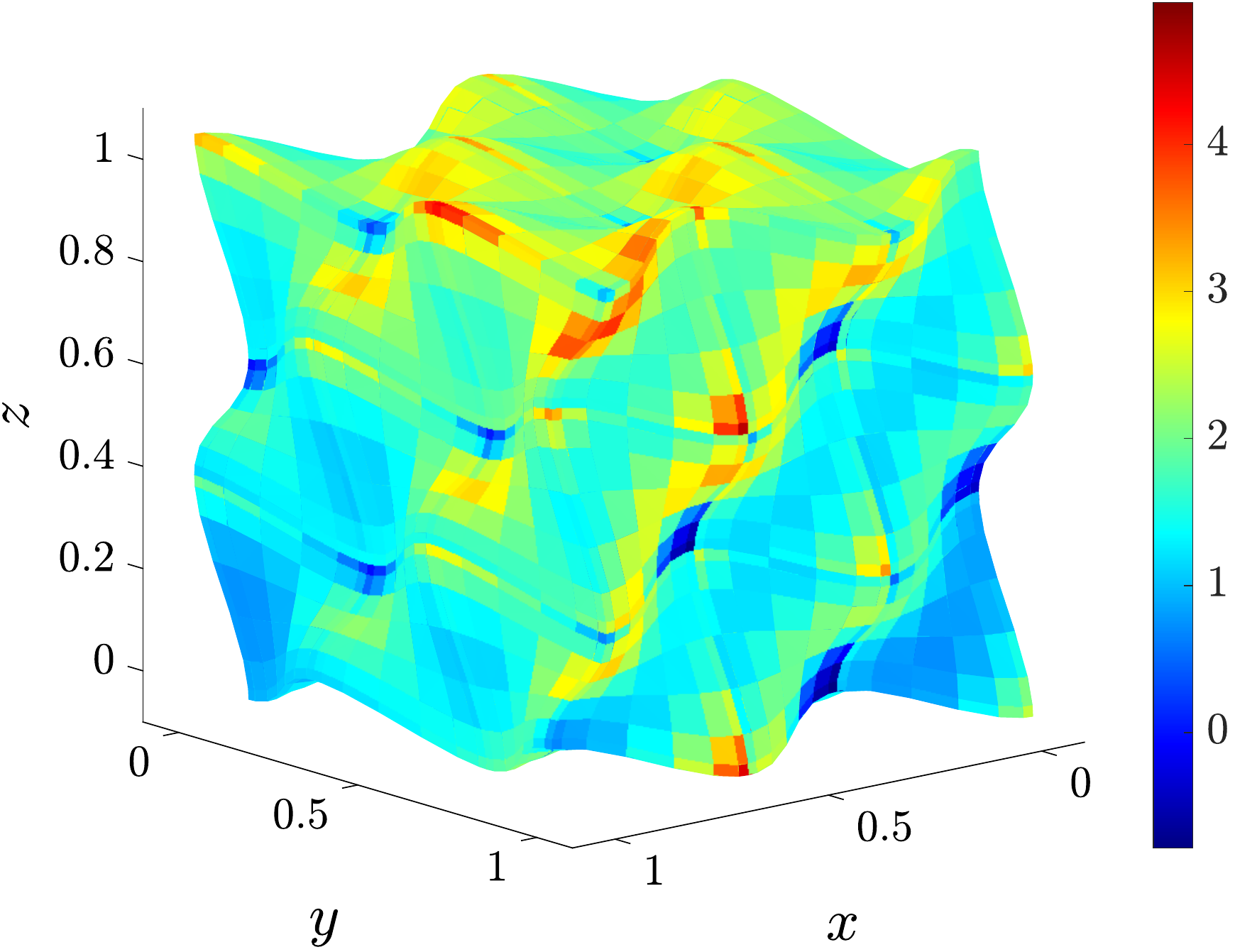}
\includegraphics[width=0.33\linewidth]{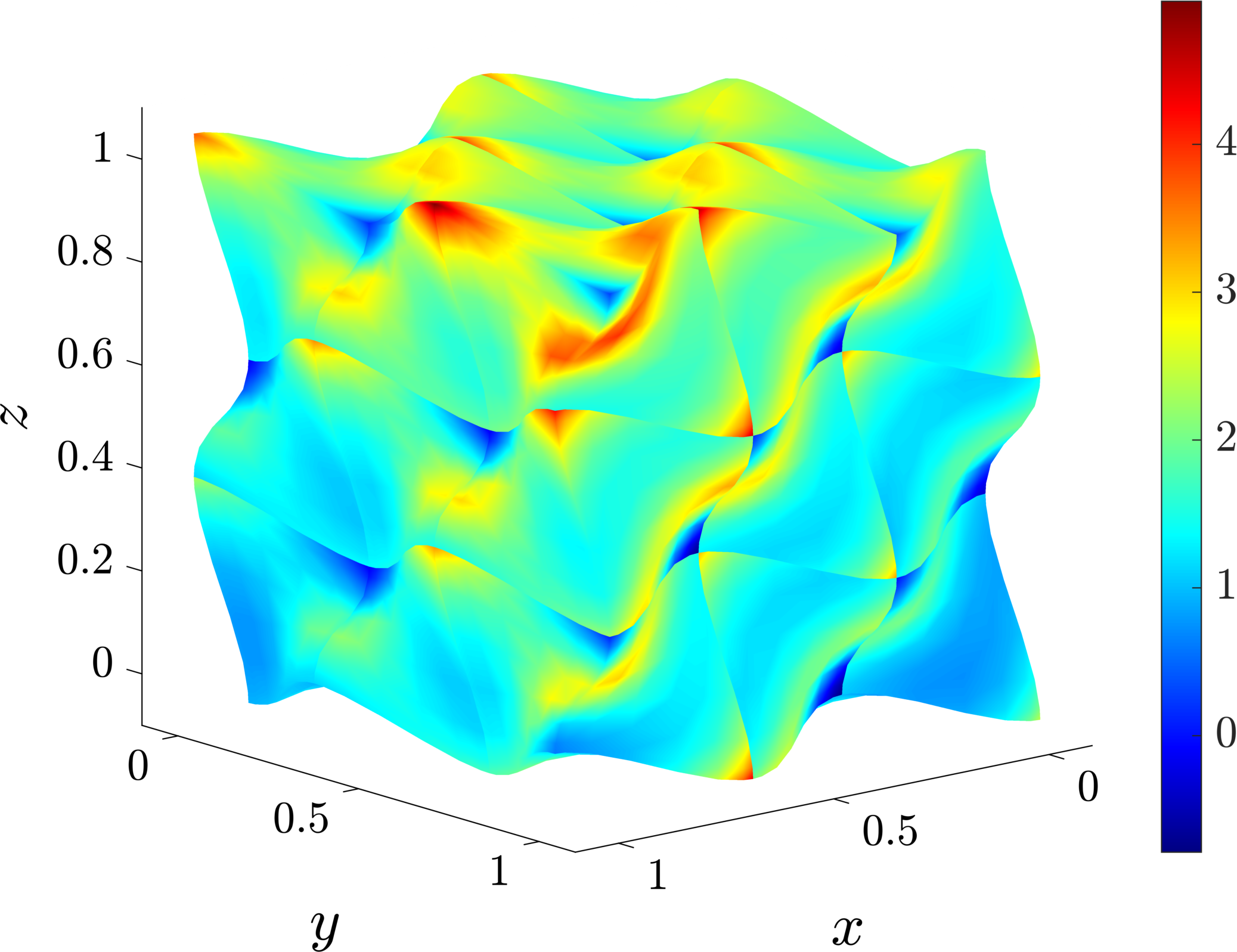}
\includegraphics[width=0.33\linewidth]{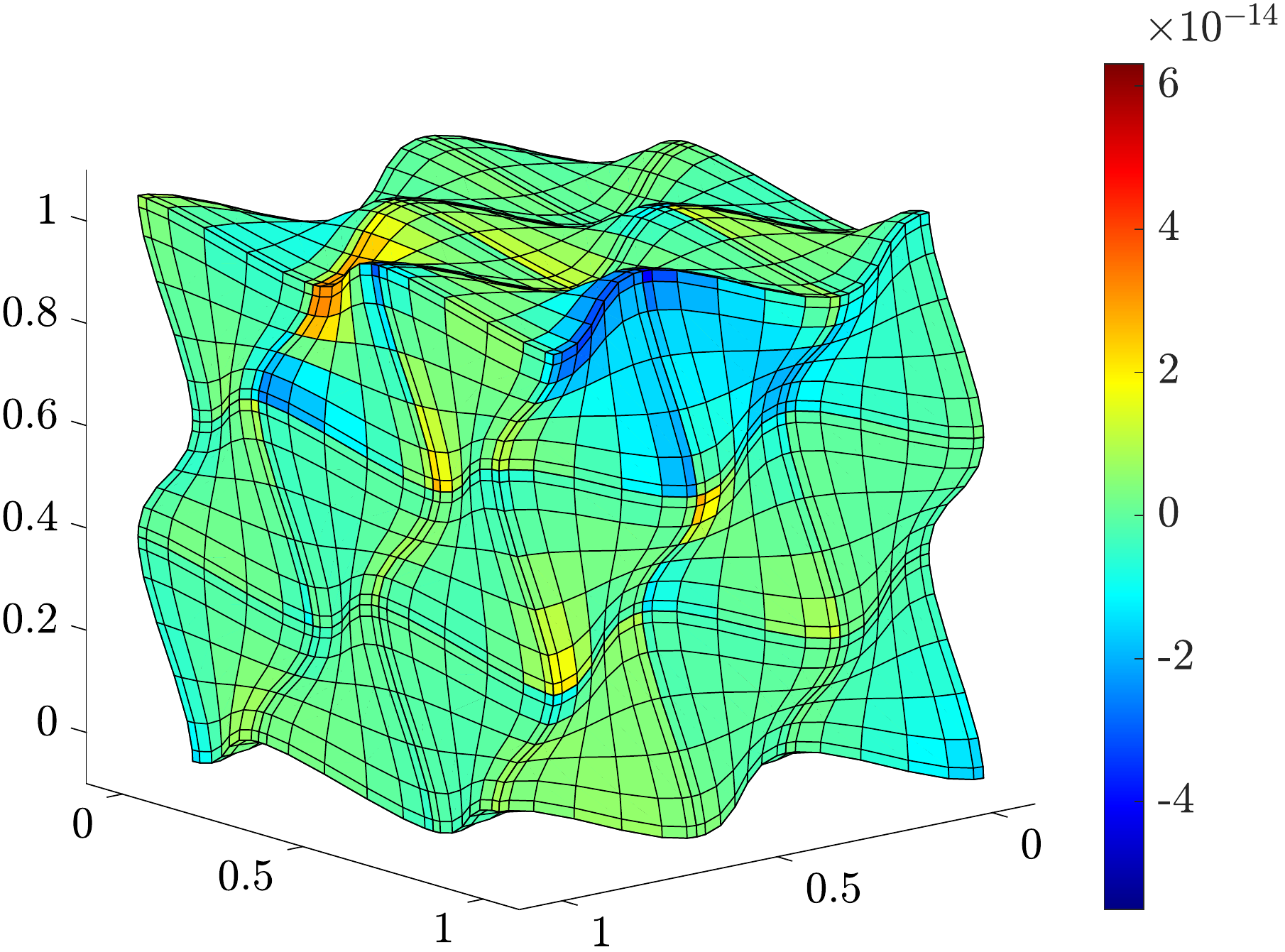} \nl
\includegraphics[width=0.33\linewidth]{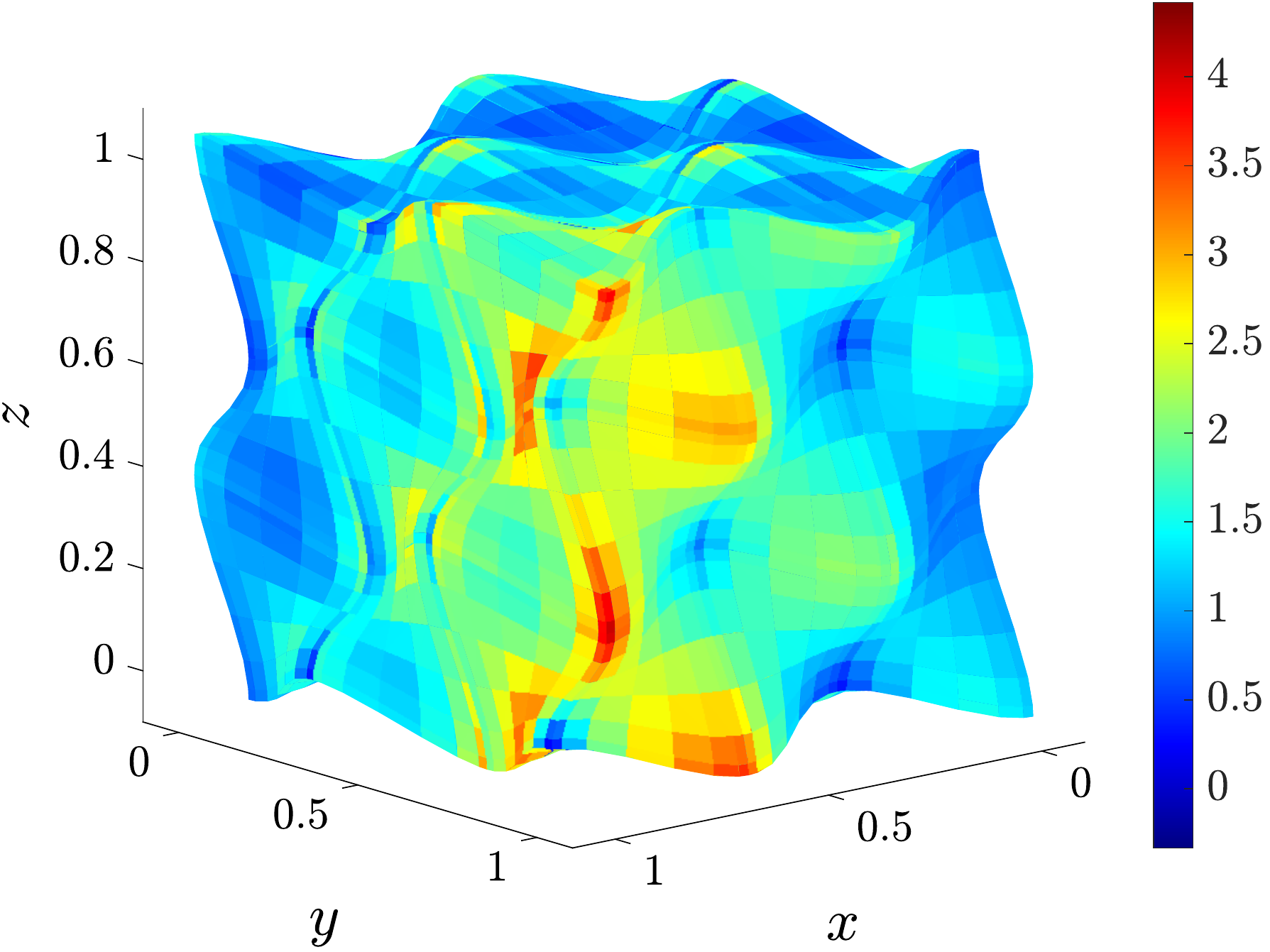}
\includegraphics[width=0.33\linewidth]{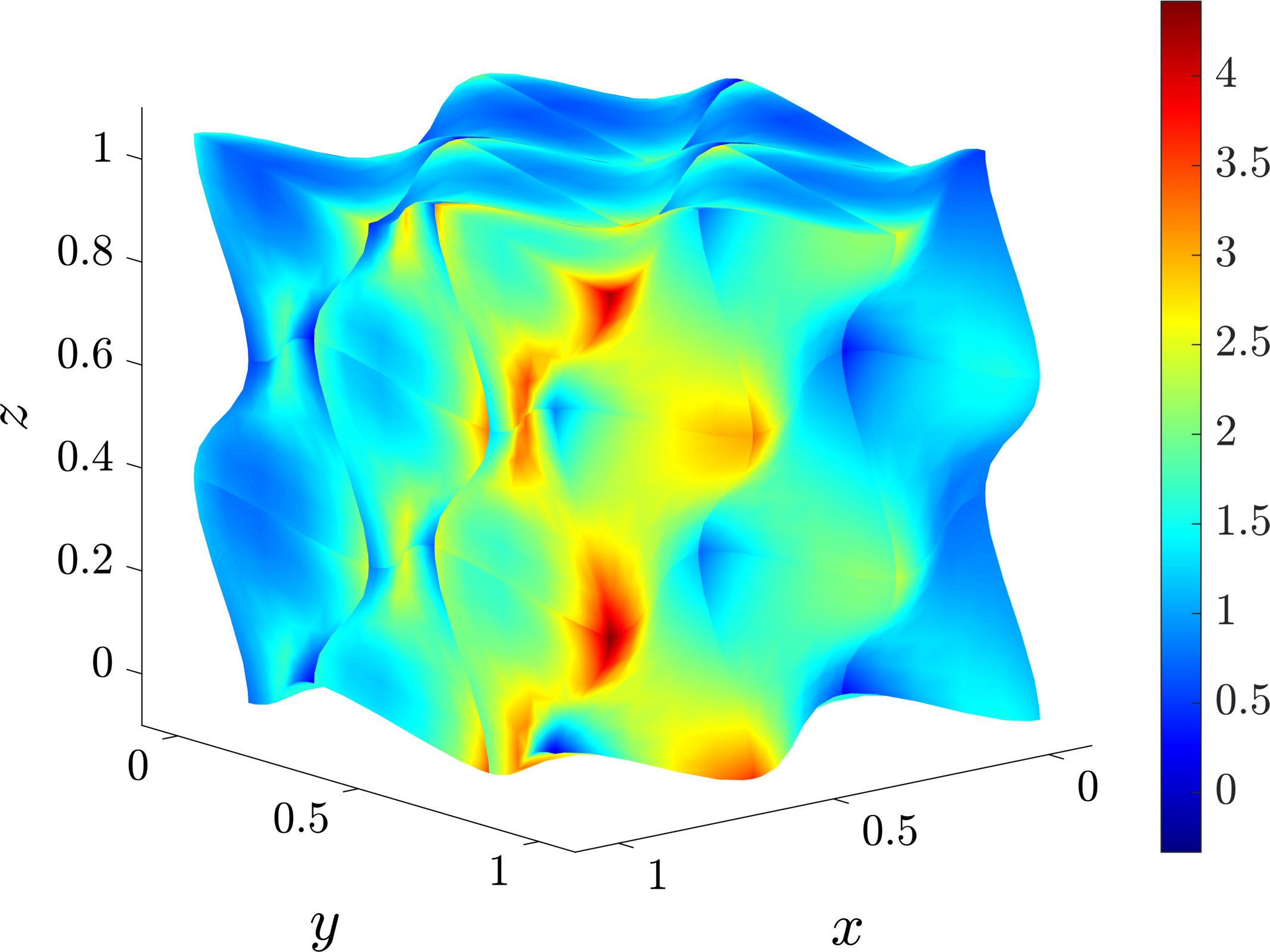}
\includegraphics[width=0.33\linewidth]{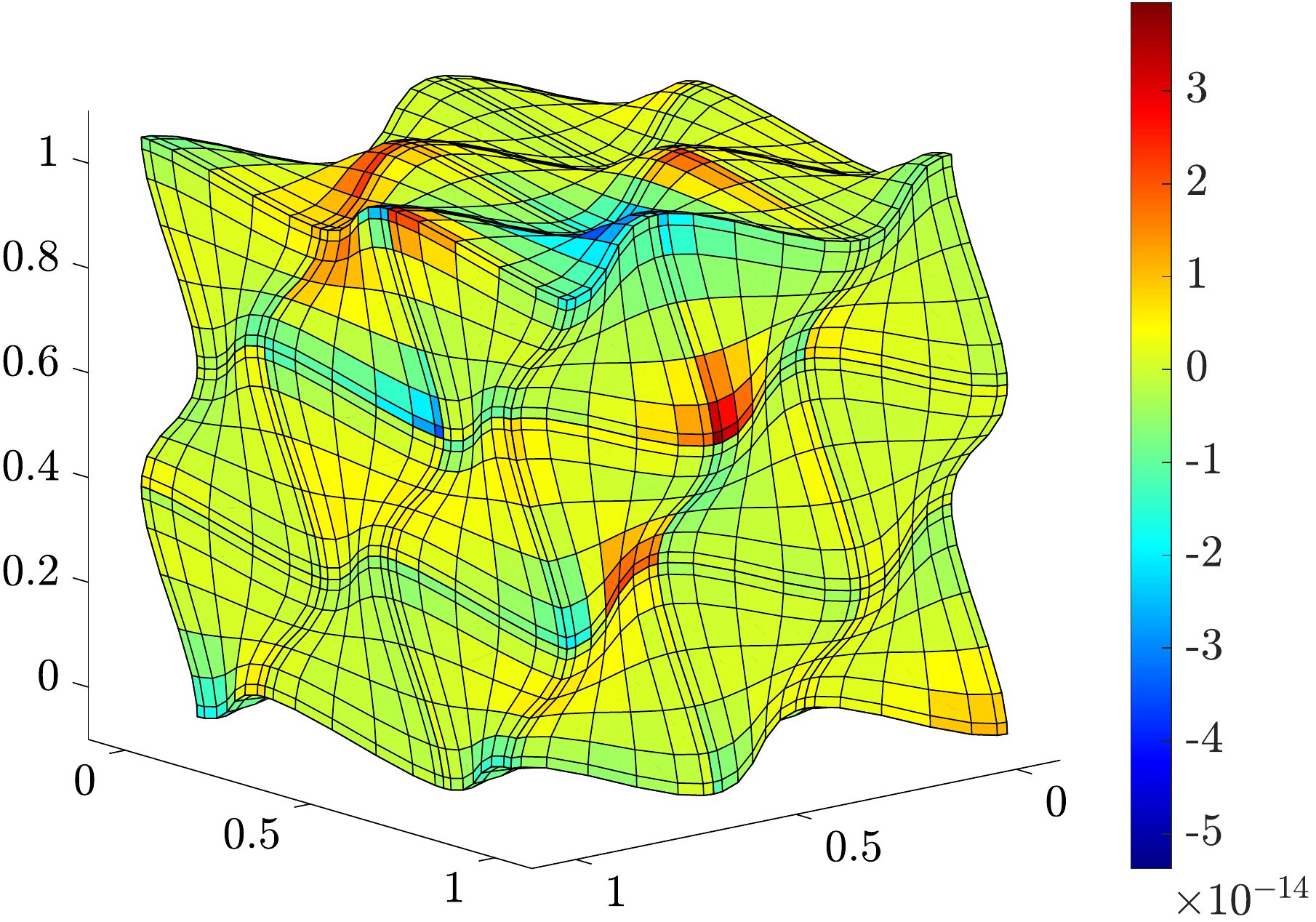}
\caption{\VJ{Comparison of pressure and velocity profiles for the continuous formulation and the DD formulation for $K=3, N=2$ case.
In the first column we see the results from continuous element formulation, in the second column we see the results from DD formulation, and in the third column we present the difference in results between both the formulations.
In the first row we see the pressure, in the second row the $x$-component of the velocity field, in the third row the $y$-component of the velocity field and in the fourth row the $z$-component of the velocity field.}}
\label{fig:comparison}
\end{figure}
\VJ{In \figref{fig:comparison} we compare the pressure and velocity profiles of the continuous formulation and the DD formulation for $K=3, N=2$ case.
In the first row we plot the pressure, in the second row we plot the $x$-component of the velocity, in the third row we plot the $y$-component of the velocity, and in the fourth row we plot the $z$-component of the velocity.
In the first column we plot the results from continuous element formulation, in the second column we plot the results from DD formulation and in the third column we plot the difference of results between both the formulations.
In the third column we see that  the maximum difference between the continuous and the DD formulation for any of the pressure or velocity profiles is $10^{-13}$.}

\begin{figure}[!htb]
	\centering
        \begin{subfigure}[b]{0.45\textwidth}
                \includegraphics[width=\textwidth]{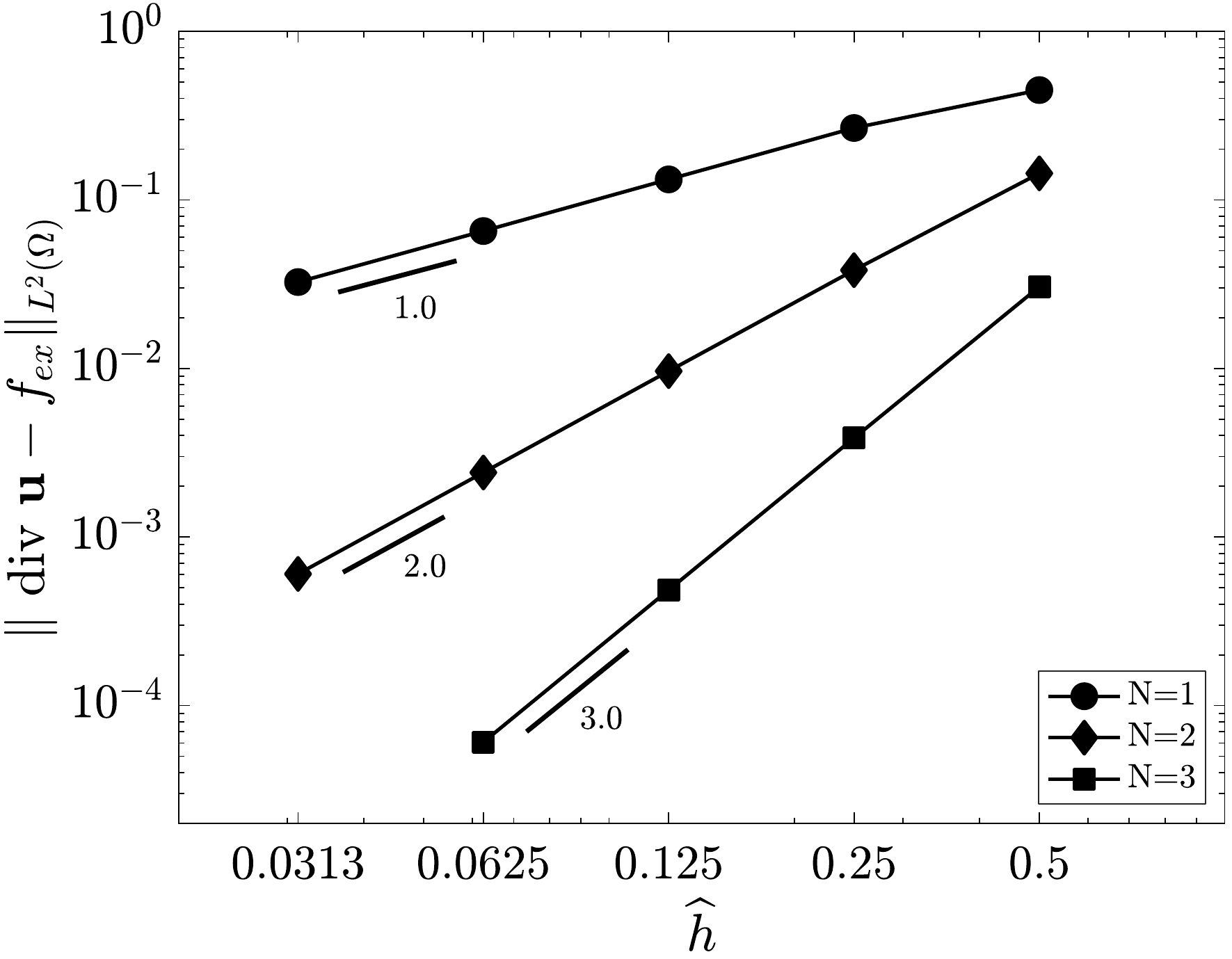}
        \end{subfigure}%
        \quad
        ~ %add desired spacing between images, e. g. ~, \quad, \qquad etc.
        \begin{subfigure}[b]{0.45\textwidth}
                \includegraphics[width=\textwidth]{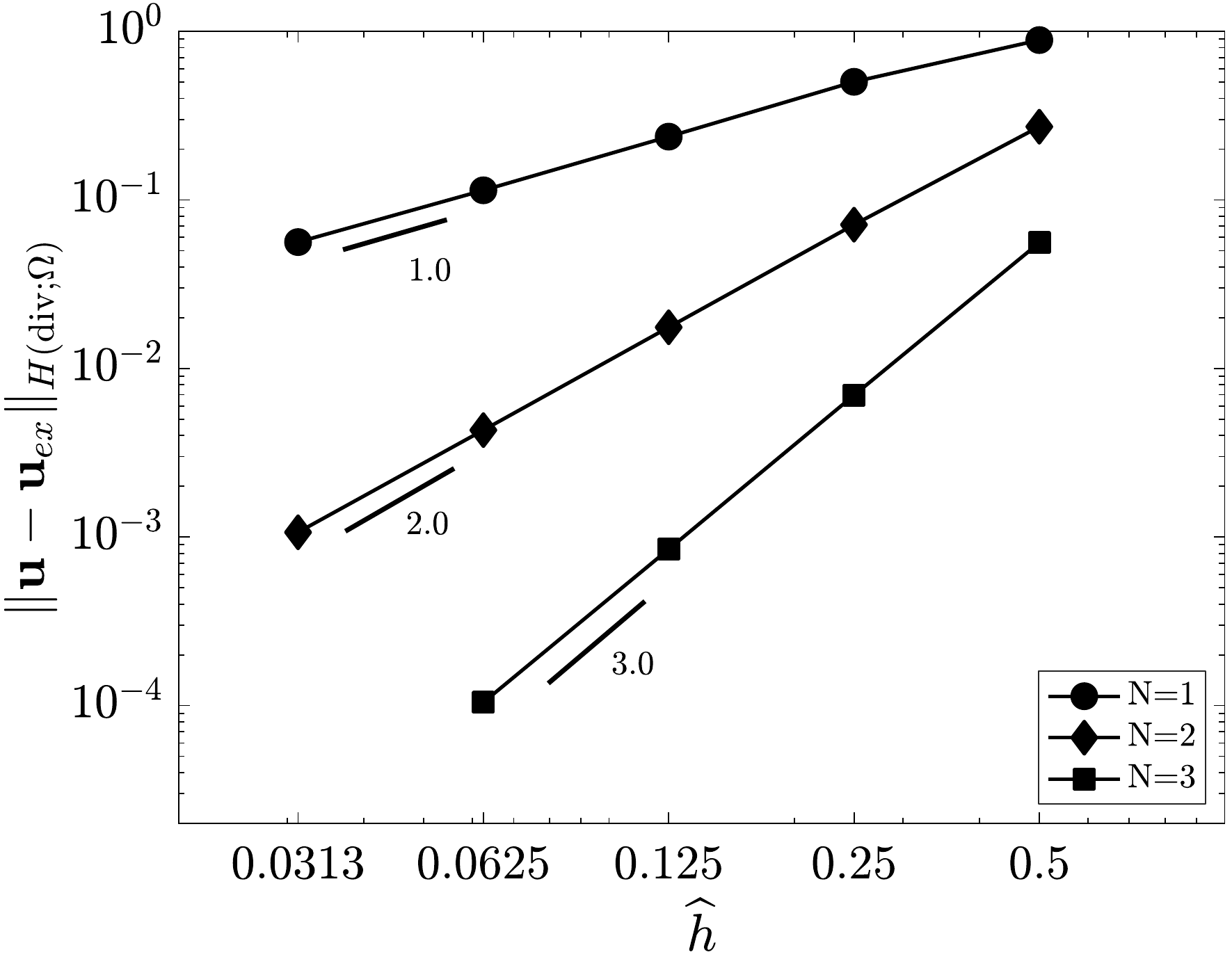}
        \end{subfigure}
        \nl
        ~ %add desired spacing between images, e. g. ~, \quad, \qquad etc.
        \begin{subfigure}[b]{0.45\textwidth}
                \includegraphics[width=\textwidth]{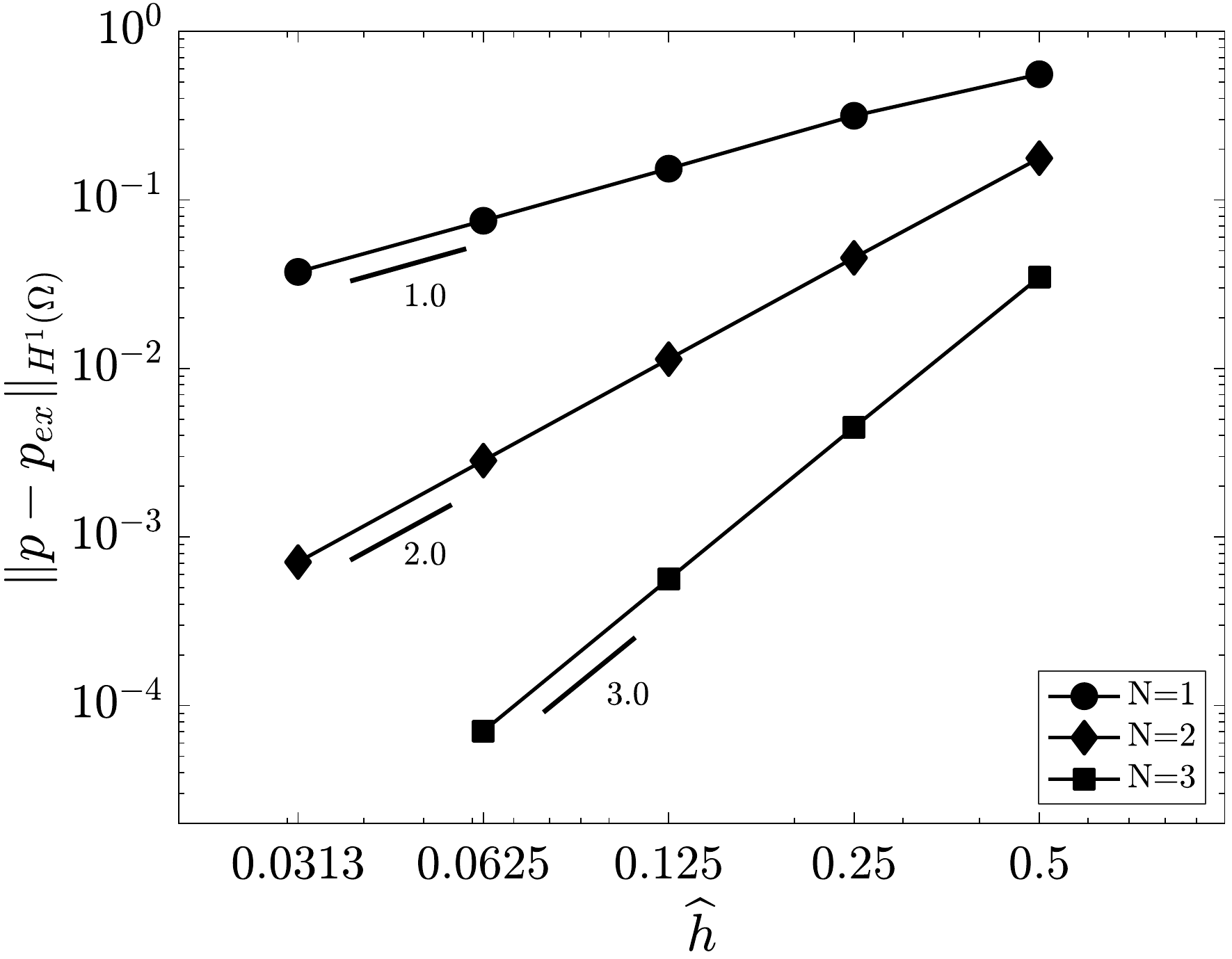}
        \end{subfigure}
		\caption{\VJ{Error convergence results for the DD formulation}. Top left: Error $(\mathrm{div}\, \bm{u} - f_{ex})$ in $L^2\bb{\Omega}$ norm. Top right: Error $(\bm{u} - \bm{u}_{ex})$ in $H\bb{\mathrm{div};\Omega}$ norm. Bottom centre: Error $(p-p_{ex})$ in $H^1\bb{\Omega}$ norm.}
		\label{fig:convergence_results_wheeler}
\end{figure}
In \figref{fig:convergence_results_wheeler} \VJ{we show the error convergence results from the DD formulation}.
At the top-left we show the convergence of the $L^2$-error for the constraint $(\mathrm{div}\ \bm{u} -f_{ex})$,  at top-right we show the convergence of the error in the $H\bb{\mathrm{div};\Omega}$ norm for the velocity field, and in the bottom-centre we show the convergence of the error in the $H^1\bb{\Omega}$ norm (using \eqref{eq:H1_norm}), for the pressure field.
On the $x$-axis we have the length of the non-deformed element $\widehat{h}$.
All the error plots show optimal rate of convergence of order $\mathcal{O}\bb{N}$.

\begin{table}[!hbt]
\centering
\caption{Average computational time (in seconds) for set-up and solution of continuous elements formulation.}
\begin{tabular}{|cc|c|c|cc|}
\hline
%& $K$   & Continuous	& \multicolumn{2}{c}{DD} &  \\
& $K$   & Set-up time	& Solve time & Total time & \\
\hline
$N=1$  & & & & & \\
\hline
& 4 		&		$<0.1$		& 	$<0.1$		& $<0.2$ 		& \\
& 8 		&		$<0.1$		& 	$<0.1$		& $<0.2$ 		& \\
& 16 		& 		0.3			& 6.5 			& 6.8 		& \\
& 32 	&		3.3			& 1327.1		& 1330.4 	& \\
& 64 	& Out of memory & & & \\
\hline
$N=2$ & & & & & \\
\hline
& 4 		& 	$<0.1$	& $<0.1$ 		& $<0.2$ 		& \\
& 8 		& 	$<0.1$	& 4.9			& 5.0 		& \\
& 16 		& 	1.3	& 1117.3 	& 1118.6	& \\
& 32 	& Out of memory & & & \\
\hline
$N=3$ & & & & & \\
\hline
& 4 		&  $<0.1$		& 0.8		& 0.9  	& \\
& 8 		& 	0.5		& 96.1		& 96.6 	& \\
& 16 	& Out of memory & & & \\
\hline
\end{tabular}
\label{tab:hybrid0}
\end{table}

\begin{table}[!hbt]
\centering
\caption{Average computational time (in seconds) for set-up and solution of DD formulation.}
\begin{tabular}{|cc|cc|c|c|c|c|c|c|}
\hline
%& $K$   & Continuous	& \multicolumn{2}{c}{DD} &  \\
& $K$ & $K1$ & $K2$  & Set-up time	& Solve \eqref{eq:lambda} & Solve \eqref{eq:pressure} \& \eqref{eq:velocity} & Total & \% to solve \eqref{eq:lambda} \\
\hline
$N=1$  & & & & & & & & \\
\hline
& 4 		& 2	& 2 	&		0.1			& $<0.1$		& $<0.1$ 		& $<0.3$ 	& - \\
& 8 		& 4 	& 2	&		0.1			& $<0.1$		& $<0.1$		& $<0.3$ 	& - \\
& 16 		& 8 	& 2	&		0.8			& $<0.1$		& $<0.1$ 		& $<1.0$ 	& - \\
& 32 	& 16 	& 2	& 		6.4			& 0.4 			& 0.1 			& 6.9 		& 5.8 \% \\
& 64 	& 32 & 2	&		50.1			& 21.3			& 0.5 			& 71.9 		& 29.6 \% \\
& 128 	& 64 & 2	&		410.3		& 2149.2		& 19.7 			& 2579.2 	& 83.3 \% \\
\hline
$N=2$ & & & & & & & & \\
\hline
& 4 		& 2	& 2 		& 	0.2		& $<0.1$ 		& $<0.1$ 		& $<0.4$ 	& - \\
& 8 		& 4	& 2 		& 	0.2		& $<0.1$		& $<0.1$ 		& $<0.4$ 	& - \\
& 16 		& 8	& 2 		& 	1.9		& 0.4			& 0.1 			& 2.3 		& 13.0 \% \\
& 32 	& 16	& 2 		& 	15.5		& 11.6 			& 0.7			& 27.9 		& 41.9 \% \\
& 64 	& 32 & 2 		& 160.9 	&  1763.5		& 63.0			& 1987.4 	& 88.7 \% \\
\hline
$N=3$ & & & & & & & & \\
\hline
& 4 		& 4	& 1 		&  0.2		& $<0.1$	& $<0.1$  	& $<0.4$ 	& - \\
& 8 		& 8	& 1 		& 	0.7		& 0.1			& $<0.1$ 	& $<0.9$ 	&  11.1 \% \\
& 16 		& 16	& 1 		& 	5.3		& 2.7			& 0.2 		& 8.2 		&  32.9 \% \\
& 32 	& 32	& 1 		& 51.0		&  267.0	& 	9.5		& 327.5 	&  81.5 \% \\
\hline
\end{tabular}
\label{tab:hybrid1}
\end{table}

\begin{table}[!hbt]
\centering
\caption{Average computational time (in seconds) for i) continuous formulation, ii) DD formulation, and the speed-up in simulation time \eqref{eq:speed_up}.}
\begin{tabular}{|cc|c|c|cc|}
\hline
%& $K$   & Continuous	& \multicolumn{2}{c}{DD} &  \\
& $K$   & Continuous	& DD & Speed-up &  \\
\hline
$N=1$  & & & & & \\
\hline
& 4 		&	$<0.2$					& 	$<0.3$		& 	-	 	&	\\
& 8 		&	$<0.2$					& 	$<0.3$		& 	-			& \\
& 16 		& 	6.8					& $<1.0$ 		&	6.8x 			&	\\
& 32 	&	1330.4				& 6.9			&	192.8x 	& \\
& 64 	&	Out of memory	& 71.9			& 	-				&	\\
\hline
$N=2$ & & & & & \\
\hline
& 4 		& 	$<0.2$					& $<0.4$ 	&  - 		& \\
& 8 		& 	5.0					& $<0.4$ 	& 12.5x			& \\
& 16 		& 	1118.6				& 2.3 		& 486.4x  	& \\
& 32 	& 	Out of memory	& 27.9 		& 	- 				& \\
\hline
$N=3$ & & & & & \\
\hline
& 4 		&  0.9					& $<0.4$	& 	-	 	& \\
& 8 		& 96.6					& $<0.9$	&	107.3x	 	& \\
& 16		& Out of memory	& 8.2		&	-				& \\
\hline
\end{tabular}
\label{tab:time_new}
\end{table}

\VJ{In \tabref{tab:hybrid0}, \tabref{tab:hybrid1}, \tabref{tab:time_new} we give the average simulation times to solve the continuous formulation, the DD formulation and the comparison between the two cases respectively.
The simulation time (in seconds) is the average for five simulation runs.
}

\[\| \bm{u} \|_{H\bb{div};\Omega}^2 = \| \bm{u} \| ^2_{L^2} + \| \mathrm{div} \bm{u} \|^2_{L^2} \;. \]

\[ \| \bm{u} - \bm{u}_{ex} \|^2 = \int_{\Omega} (\bm{u} - \bm{u}_{ex})^2 \mathrm{d}\Omega \]

\VJ{In \tabref{tab:hybrid0}, we give the average simulation time to solve for continuous formulation for elements of order $N=1,2,3$.
In the first column we have the number of elements in one direction, in the second column we have the time taken to set-up the matrices and in the third column the time taken to solve \eqref{eq:cont_p} and \eqref{eq:cont_u}.
The simulation times that are less than $0.1$ second, have not been measured with further accuracy.
For $N=1$ the system run out of memory for $K=64$, for $N=2$ at $K=32$, and for $N=3$ at $K=16$.}

\VJ{In \tabref{tab:hybrid1}, we give the average simulation time to solve for DD formulation for elements of order $N=1, 2,3$.
In the first column, we give the total numer of elements in one direction.
In the second column we give the number of sub domains in one direction.
In the third column we give the number of elements in each direction within each sub domain.
In the fourth column we give the time to set-up the matrices.
In the fifth column we give the time to solve \eqref{eq:lambda}.
In the sixth column we give the time to solve \eqref{eq:velocity} and \eqref{eq:pressure}.
In the seventh column we give the total time, and in the eigth column we give the \% amount of time spent to evaluate \eqref{eq:lambda}.
For the cases where simulation times are less than 0.1 second the accurate simulation times are not determined and therefore we do not calculate the last column.
We observe in the last column that for all the three cases i.e. $N=1,2,3$, as we increase the total number of elements the percentage of time spent on solution of the Lagrange multiplier system increases.
In the last case listed in the table for $N=1,2,3$, the time spent in solution of \eqref{eq:lambda} is above $80 \%$.
Given that DD formulation is for large simulations, we see that in this case the majority of the time is spent on solution of the Lagrange multiplier system \eqref{eq:lambda}.
}

\VJ{In \tabref{tab:time_new} we compare the total time taken to solve the continuous formulation and the DD formulation for elements of order $N=1,2,3$.
The comparison is made of cases with same total number of elements in one direction, such that topology of the mesh discretization remains the same.}
In the first column we have the total number of elements in one direction $K$, in the second column we have the total simulation time from continuous elements formulation, in the third column we have the total simulation time from DD formulation.
In the fourth column we give the speed-up factor as
\begin{equation} \label{eq:speed_up}
\mbox{speed-up}  = \frac{simulation\ time\ continuous\ formulation}{simulation\ time\ DD\ formulation}  \;.
\end{equation}
%We present the results for elements of order $N=1,2,3$.
Firstly, we see that for all $N$, continuous formulation runs out of memory for lower values of $K$.
That is, the DD formulation requires relatively less memory than the continuous formulation.
%the memory requirement for DD method is less than that for continuous formulation.
%In case of DD methods we can go at least 2 times more refined mesh.
%This is evident because the global system of continuous formulation is much larger in size as compared to global $\lambda$ system \eqref{eq:lambda}.
Secondly, the simulation time required, for the same refinement, is less for DD method than for continuous formulation in all cases.
Thirdly, the speed-up factor increases as the mesh is more refined, i.e. the total number of elements are increased.
%In the last column on the right, we show the gain in time we achieve using DD method.
%We observe that the gains are higher for high order methods, for example for N=3, K=4 case the DD simulation is 24949 times faster, which is because high order elements have more internal degrees of freedom and less boundary degrees of freedom.

\begin{figure}
\centering
\includegraphics[scale=0.5]{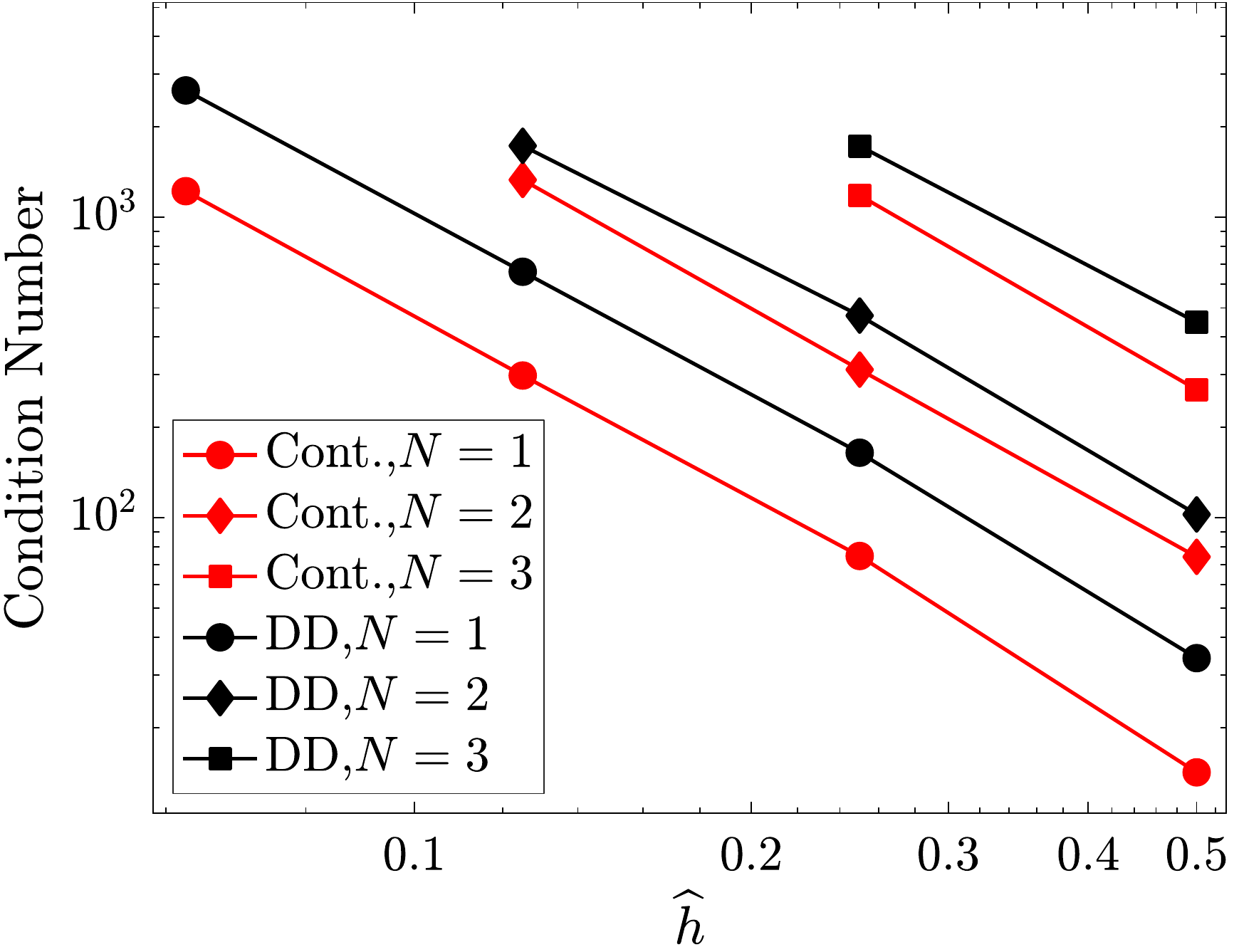}
\caption{Comparison of condition number of \eqref{eq:cont_p} and \eqref{eq:lambda} for $N=1,2,3$.}
\label{fig:cond}
\end{figure}
\VJ{In \figref{fig:cond} we compare the condition number of the global system of \eqref{eq:cont_p} and \eqref{eq:lambda}.
We observe that in both the cases, for $N=1,2,3$, and same discretization, the condition number of the Lagrange multiplier system \eqref{eq:lambda} is higher than the condition number of the pressure unknowns of \eqref{eq:cont_p}.
Moreover, the rate of growth of condition number for $N=1,2,3$ is similar in both the cases.}
\subsection{Test case II: SPE 10}
\begin{figure}
\centering
\includegraphics[scale=0.4]{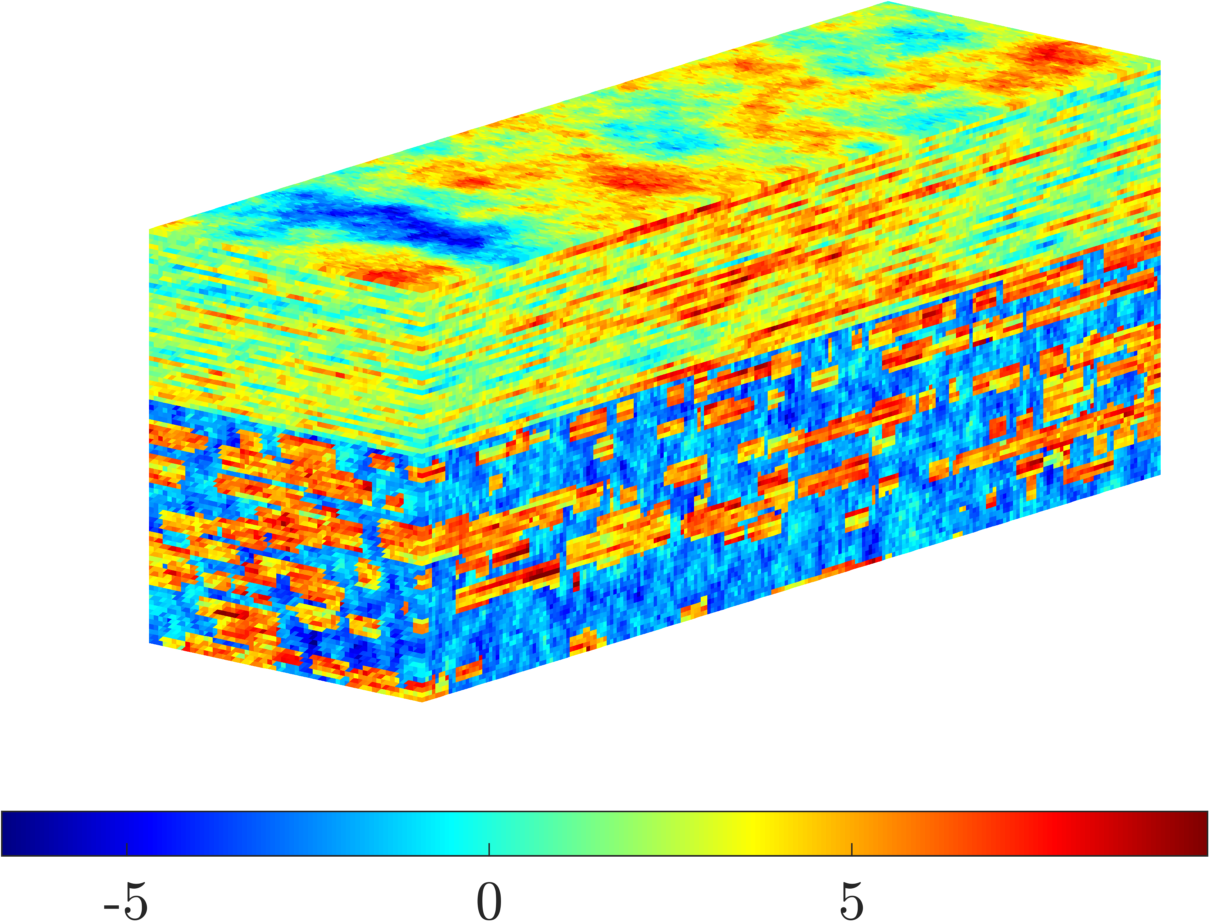}
\caption{\VJ{Natural logarithm of permeability field for SPE10 case.}}
\label{fig:SPE_domain}
\end{figure}
In this section we show the results for solution of SPE10 benchmark problem~\cite{spe10}.
It is often used for validation of numerical schemes for reservoir modelling applications because of its challenging permeability field.
The domain of the problem and the natural logarithm of the permeability field is shown in \figref{fig:SPE_domain}.
The size of the domain is $1200ft \times 2200ft \times 170ft$.
The domain is divided into equal blocks of size $20ft \times 10ft \times 2ft$, i.e. $60\times 220 \times 85 = 1122000$ blocks.
Each block has a constant isotropic permeability tensor \VJ{$\mathbb{K}_i = k_i\mathbb{I}$}.
%The contours in \figref{fig:SPE_domain} show the natural logarithm of the permeability field component in \VJ{all} direction.
The top $70ft$ of the domain represents the Tarbert formation that has smooth changes in the permeability field component, and the bottom $100ft$ represents the Upper Ness formation that has sharp changes in the permeability field component.

The right hand side term is $f_{ex} = 0$.
We impose Neumann boundary conditions, $\hat{u} = \bm{u} \cdot \bm{n} = 0$, on \VJ{$\hat{y}, \hat{z} = 0,1$}, and the Dirichlet boundary conditions, $\hat{p} = 1$ on \VJ{$\hat{x}=0$}, and $\hat{p} = 0$ on \VJ{$\hat{x}=1$}.

The numerical solution of this problem using the continuous formulation was not possible because the system ran out of memory and therefore we only present the results with DD method.
This shows that DD method can also \VJ{be advantageous in cases to reduce the memory requirements for large simulations.}

\begin{figure}
\centering
\includegraphics[scale=0.135]{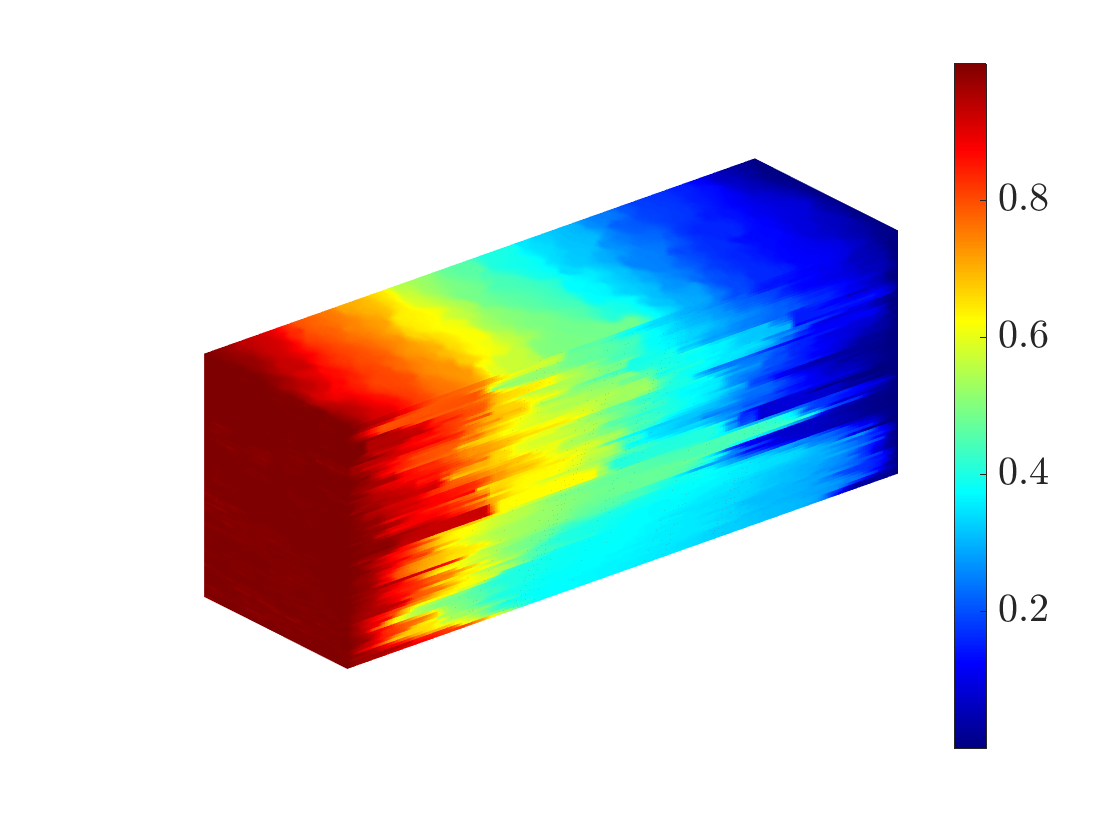} 
\includegraphics[scale=0.135]{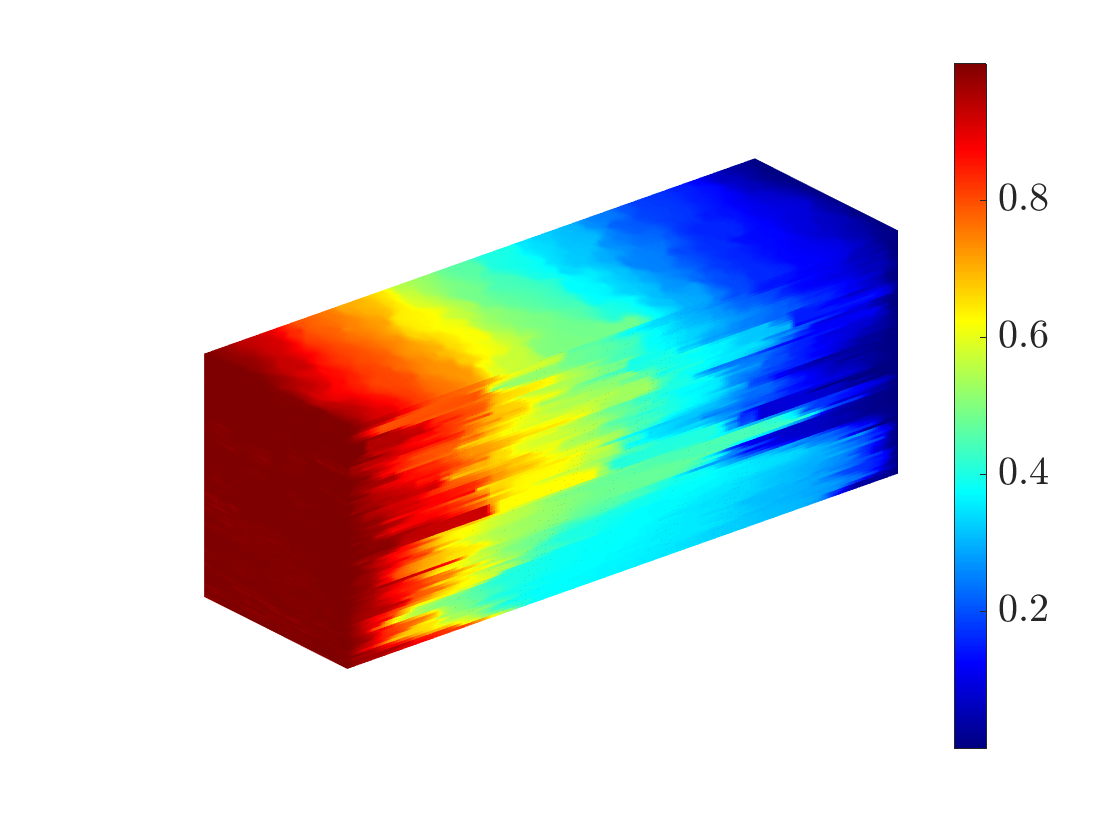} 
\includegraphics[scale=0.135]{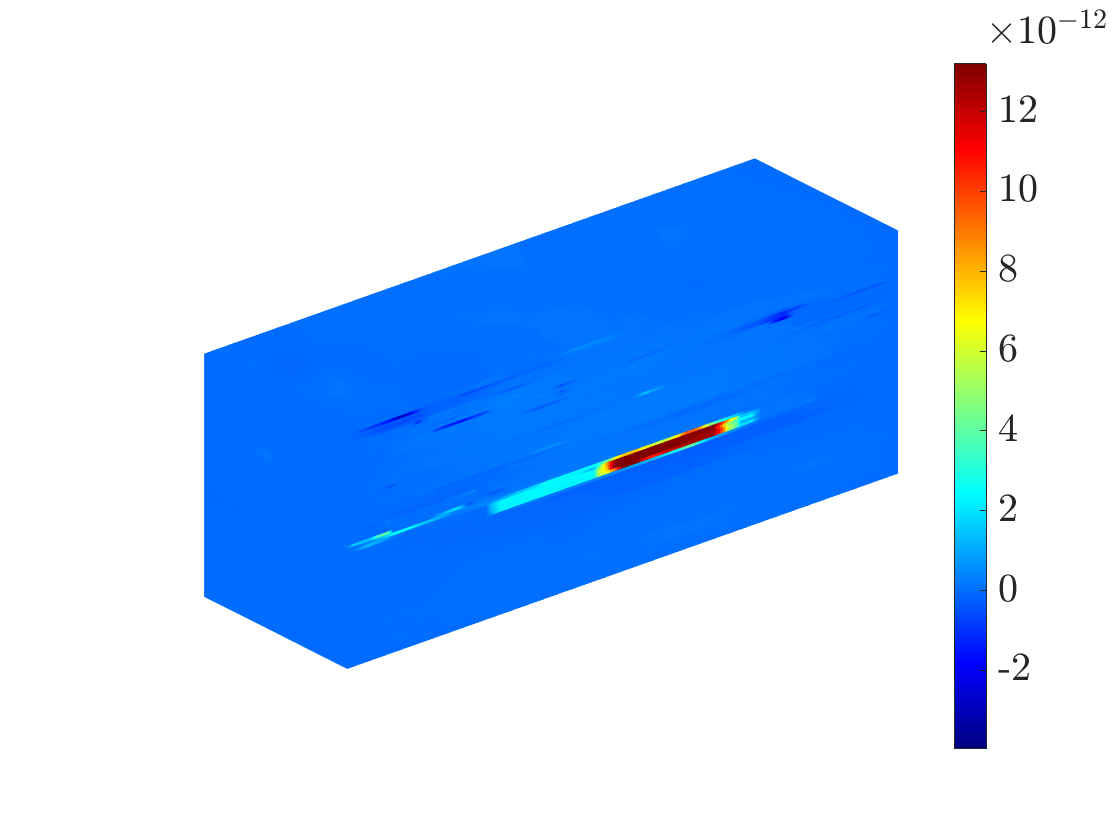}
\caption{The contour plots of pressure field for SPE10 benchmark problem. Left: results from $Case\ 1$. Centre: results from $Case\ 2$. Right: difference between results from $Case\ 1$ and $Case\ 2$.}
\label{fig:result_pressure}
\end{figure}
\begin{table}[!hbt]
\centering
\caption{Average computational time (in seconds) for set-up and solution of SPE10 case using DD formulation.}
\begin{tabular}{|c|c|cc|c|c|c|c|c|c|}
\hline
%& $K$   & Continuous	& \multicolumn{2}{c}{DD} &  \\
& $K$ & $K1$ & $K2$  & Set-up time	& Solve \eqref{eq:lambda} & Solve \eqref{eq:pressure} \& \eqref{eq:velocity} & Total & \% \eqref{eq:lambda} \\
\hline
Case  1 & 1122000 	& $15 \times 55 \times 17 $	& $4 \times 4 \times 5 $ 	&	154.0	&  259.4	& 20.4 		& 433.8 	& 59.8 \% \\
Case 2 & 1122000 	& $12 \times 44 \times 17$ 	& $5 \times 5 \times 5 $	&	184.2 	& 159.5		& 24.3		&  368.0	& 43.3 \% \\
\hline
\end{tabular}
\label{tab:spe10}
\end{table}
For numerical solution using the DD method, we use two different approaches.
In $Case\ 1$, we divide the domain into $15 \times 55 \times 17$ sub domains, and each sub domain is further divided into $4 \times 4 \times 5$ elements.
In $Case\ 2$, we divide the domain into $12 \times 44 \times 17$ sub domains, and each sub domain is further divided into $5 \times 5 \times 5$ elements.
\VJ{We also considered a third case with $6 \times 22 \times 17$ sub domains, and each sub domain divided into $10 \times 10 \times 5$ elements, but this case exceeded the memory bounds of the system.}
The total number and topological configuration of the elements remain the same in all the cases.
We use only the lowest order elements, i.e. $N=1$, in both the cases.
The results of the pressure field from both the cases are shown in \figref{fig:result_pressure}.
In the left plot we see the results from $Case\ 1$, in the centre plot we see the results from $Case\ 2$, and in the right plot we see the difference between the two results.
We observe that the maximum difference in results is of order $10^{-12}$.

\VJ{In \tabref{tab:spe10} we give the average run time of five simulation runs to solve the  SPE10 case using DD formulation.}
In the second column we have the total number of elements, in the third column we have the total number of sub domains, in the fourth column we have the total number of elements within each sub domain, in the fifth column we have the set-up time of matrices, in the sixth column we have the time taken to solve \eqref{eq:lambda}, in the seventh column we have the time taken to solve $p$ \& $\bm{u}$. \eqref{eq:pressure} \& \eqref{eq:velocity}, in the eighth column we have the total time, and in the ninth column we have the \% of time taken to solve the Lagrange multiplier system \eqref{eq:lambda} with respect to the total time.
%The average run time, of five simulation runs, to solve the pressure and velocity unknowns using DD formulation in $Case 1$ is $471s$ seconds, and $Case 2$ is $406s$ seconds, rounded off to nearest integer.
%The condition number for the $\lambda$ system for $Case 1$ is xxx, and for $Case 2$ is xxx.
The results show that for the same total number and topology of elements, the choice of sub domain decomposition also affects the simulation times.
% and condition number of the $\lambda$ system.
%\VJ{!!! It also affects the memory requirements!!!}.

This test case demonstrates that use of algebraic dual representations with DD method can be applied to practical applications.
\section{Conclusions} \label{sec:conclusions}
In this paper we have presented the use of algebraic dual spaces for DD method for Darcy flow.
We have defined the broken Sobolev spaces and their finite dimensional counterparts.
%, $D\bb{\Omega _T}$, $S\bb{\Omega _T}$ and their dual counterparts for broken Sobolev spaces $H\bb{\mathrm{div};\Omega_T}$ and $L^2\bb{\Omega _T}$. 
We have also defined the global finite dimensional trace space that connects the broken spaces.
These spaces are used to solve the DD formulation of Darcy equations.
\VJ{It is shown that using algebraic dual \representations the matrix representation of continuity constraint across the sub domains becomes sparse and metric-free.}
The first test case is a manufactured solution where it is shown that i) the results from continuous formulation and DD formulation are precisely the same and that DD formulation also has optimal rate of convergence of errors, ii) the DD formulation is more memory efficient, iii) the DD formulation is more efficient in terms of simulation run times.
In the second test case we show that DD scheme can be used to reduce the memory requirements and solve for large practical applications.

We have demonstrated that algebraic dual spaces can be used with DD schemes.
In future this work will be extended to broken $H^1\bb{\Omega}$ and $H\bb{\mathrm{curl};\Omega}$ spaces, to address problems such as vector Laplacian, Stokes flow, and linear elasticity.

\section{Acknowledgements}
The authors would like to thank Dr. Hadi Hajibeygi for the processed permeability data for SPE10 benchmark case.

\bibliographystyle{elsarticle-num}
\bibliography{references}

\end{document}